\documentclass[11pt,fullpage]{article}

\usepackage{amsmath}
\usepackage{amssymb}
\usepackage{amsfonts}
\usepackage{xr}

\numberwithin{equation}{section}
\newtheorem{thm}{Theorem}[section] 
\newtheorem{prp}[thm]{Proposition}
\newtheorem{lmm}[thm]{Lemma}   
\newtheorem{crl}[thm]{Corollary} 
\newtheorem{dfn}[thm]{Definition} 

\def\e_ref#1{(\ref{#1})}
\def\under#1{\underline{#1}}
\def\ov#1{\overline{#1}}
\def\lra{\longrightarrow}
\def\Lra{\Longrightarrow}

\def\lan{\langle}
\def\ran{\rangle}

\def\al{\alpha}
\def\de{\delta}

\def\ga{\gamma}
\def\io{\iota}
\def\ka{\kappa}
\def\la{\lambda}
\def\om{\omega}
\def\si{\sigma}
\def\th{\theta}
\def\ve{\varepsilon}
\def\ups{\upsilon}

\def\vt{\vartheta}

\def\De{\Delta}
\def\Ga{\Gamma}
\def\Om{\Omega}
\def\Si{\Sigma}
\def\Th{\Theta}

\def\codim{\text{codim}}
\def\rk{\text{rk}}
\def\ev{\text{ev}}
\def\RT{\text{RT}}
\def\Im{\text{Im}}

\def\P{\Bbb{P}^n}
\def\PP{\Bbb{P}^2}
\def\PPP{\Bbb{P}^3}
\def\i{\infty}
\def\eset{\emptyset}

\begin{document}

\title{Enumeration of One-Nodal Rational Curves\\
in Projective Spaces }

\author{Aleksey Zinger
\thanks{Partially supported by NSF grant DMS-9803166}}
\date{\today}

\maketitle

\begin{abstract}
\noindent
We give a formula computing the number of one-nodal rational curves 
that pass through an appropriate collection of constraints in
a complex projective space.
We combine the methods and results from three different papers.
\end{abstract}

\tableofcontents

\section{\bf Introduction} 

\noindent
Enumerative algebraic geometry is a field of mathematics that dates back
to the nineteenth century.
However, many of its most fundamental problems remained unsolved
until the early~1990s.
For example, let $d$ be a positive integer and 
$\mu\!=\!(\mu_1,\ldots,\mu_N)$ an $N$-tuple of linear subspaces of~$\P$
of codimension at least two such~that
$$\codim_{\Bbb{C}}\mu\equiv
\sum_{l=1}^{l=N}\codim_{\Bbb{C}}\mu_l-N
=d(n+1)+n-3.$$
If the constraints $\mu$ are in general position,
denote by $n_d(\mu)$ the number of rational degree-$d$ curves
that pass through $\mu_1,\ldots,\mu_N$.
This number is finite and depends only on the homology classes
of the constraints.
If $d\!=\!1$, it can be computed using Schubert calculus; see~\cite{GH}.
All but very-low-degree numbers $n_d(\mu)$ remained unknown until 
\cite{KM} and~\cite{RT} derived a recursive formula for these numbers.
In this paper, we prove

\begin{thm}
\label{main_thm}
Suppose $n\!\ge\!3$, $d\!\ge\!1$, 
and $\mu\!=\!(\mu_1,\ldots,\mu_N)$ is 
an $N$-tuple of proper subvarieties of~$\P$ in general position
such~that
\begin{equation}\label{main_thm_e}
\codim_{\Bbb{C}}\mu\equiv
\sum_{l=1}^{l=N}\codim_{\Bbb{C}}\mu_l-N=d(n+1)-1.
\end{equation}
Then the number of degree-$d$ rational curves that have a simple node
and pass through the constraints~$\mu$ is given~by
\begin{gather*}
n_d^{(1)}(\mu)=\frac{1}{2}\big(
\RT_{1,d}(\mu_1;\mu_2,\ldots,\mu_N)-CR_1(\mu)\big),
\qquad\hbox{where}\\
CR_1(\mu)=\sum_{k=1}^{2k\le n+1}\!\! (-1)^{k-1}(k\!-\!1)!
\!\sum_{l=0}^{n+1-2k} \binom{n\!+\!1}{l}
\big\lan a^l\eta_{n+1-2k-l},\big[\bar{\cal V}_k(\mu)\big]\big\ran.
\end{gather*}
The symplectic invariant $\RT_{1,d}(\cdot;\cdot)$
and the top intersections 
$\big\lan a^l\eta_{n+1-2k-l},\big[\bar{\cal V}_k(\mu)\big]\big\ran$
are computable via algorithms described elsewhere.
\end{thm}

\begin{center}
\begin{tabular}{||c|c|c|c|c|c||}
\hline\hline
$n$&                     3&       4&         5&          5&             6\\
\hline
$d$&                     4&       4&         4&          6&             6\\ 
\hline 
$\mu$&               (5,5)& (5,1,4)& (5,1,0,4)&  (2,1,1,7)&   (2,1,1,1,6)\\
\hline 
$n_d^{(1)}(\mu)$&  1,800&   1,800&     1,800&     20,340&          20,340\\
\hline\hline
\end{tabular}
\end{center}

\vspace{.2in}

\noindent
For the purposes of this table, we assume that the constraints 
$\mu_1,\ldots,\mu_N$ are linear subspaces of~$\P$ of codimension at
least~two.
We describe such a tuple~$\mu$ of constraints by listing 
the number of linear subspaces of codimension $2,\ldots,n$
among $\mu_1,\ldots,\mu_N$.\\

\noindent
In the statement of Theorem~\ref{main_thm},
$\RT_{1,d}(\cdot;\cdot)$ denotes the genus-one degree-$d$
symplectic invariant of~$\P$ defined in~\cite{RT}.
This invariant can be expressed in terms of the numbers~$n_d(\cdot)$;
see~\cite{RT}.
In particular, it is computable.
Brief remarks concerning the meaning of~$\RT_{1,d}(\cdot;\cdot)$
can be found at the beginning of Section~\ref{comp_sect}.\\

\noindent
The compact oriented topological manifold~$\bar{\cal V}_k(\mu)$ consists 
of unordered $k$-tuples of stable rational maps of total degree~$d$.
Each map comes with a special marked point~$\i_i$.
All these marked points are mapped to the same point in~$\P$.
In particular, there is a well-defined evaluation map
$$\ev\!: \bar{\cal V}_k(\mu)\lra\P,$$
which sends each tuple of stable maps to the value at 
one of the special marked points.
We also require that the union of the images of the maps in each tuple
intersect each of the constraints \hbox{$\mu_1,\ldots,\mu_N$}.
In fact, the elements in the tuple carry a total of~$N$
marked points, $y_1,\ldots,y_N$,  in addition to the $k$ 
special marked points.
These marked points are mapped to the constraints $\mu_1,\ldots,\mu_N$,
respectively.
Roughly speaking, each element of $\bar{\cal V}_k(\mu)$ corresponds 
to a degree-$d$ rational curve in~$\P$, which has at least $k$ irreducible
components, and $k$ of the components meet at the same point in~$\P$.
The precise definition of the spaces $\bar{\cal V}_k(\mu)$
can be found in Subsection~\ref{notation}.\\

\noindent
The cohomology classes~$a$ and~$\eta_l$ are tautological classes
in~$\bar{\cal V}_k(\mu)$.
In~fact,
$$a=\ev^*c_1\big({\cal O}(1_{\P})\big).$$
Let $\bar{\cal V}_k'(\mu)$ be the oriented topological manifold
defined as~$\bar{\cal V}_k(\mu)$, except without specifying
the marked points \hbox{$y_1,\ldots,y_N$} mapped to the constraints
\hbox{$\mu_1,\ldots,\mu_N$}.
Then, there is well-defined forgetful~map,
$$\pi\!:\bar{\cal V}_k(\mu)\lra\bar{\cal V}_k'(\mu),$$
which drops the marked points $y_1,\ldots,y_N$
and contracts the unstable components.
The cohomology class $\eta_l\!\in\!H^{2l}(\bar{\cal V}_k(\mu))$
is the sum of all degree-$l$ monomials in the elements of 
the set
$$\big\{\pi^*\psi_{\hat{0}_1},\ldots,\pi^*\psi_{\hat{0}_k}\big\}
\subset H^2(\bar{\cal V}_k(\mu)).$$
As common in algebraic geometry, 
$\psi_{\hat{0}_i}$ denotes the first chern class
of the universal cotangent line bundle for 
the marked point~$\hat{0}_i\!\in\!\bar{\cal V}_k'(\mu)$.
In Subsection~\ref{notation}, we give a definition of~$\eta_l$
that does not involve the projection map~$\pi$.
An algorithm for computing the intersection numbers 
involved in the statement of Theorem~\ref{main_thm} is given 
in Subsection~5.7 of~\cite{Z2}.
It is closely related to the algorithm of~\cite{P2}
for computing intersections of tautological classes in
moduli spaces of stable rational maps into~$\P$.\\

\noindent
If $n\!=\!2$, we denote by $n_d^{(1)}(\mu)$ the number of
rational degree-$d$ curves passing through
the constraints counted with a choice of the node on each curve.
The formula of Theorem~\ref{main_thm} gives
\begin{equation}
\label{n2nums}
n_d^{(1)}(\mu)=\binom{d\!-\!1}{2}n_d(\mu).
\end{equation}
This identity is clear, since the arithmetic genus
of every degree-$d$ curve in $\PP$ is~$\binom{d-1}{2}$.
Equation~\e_ref{n2nums} is used in~\cite{P1} to count
genus-one plane curves with complex structure fixed.
More precisely, if $\mu$ is a tuple of constraints in~$\P$
satisfying condition~\e_ref{main_thm_e},
let $n_{1,d}(\mu)$ denote the number of genus-one degree-$d$ curves
that pass through the constraints~$\mu$ and have 
a fixed generic complex structure on the normalization,
i.e.~its $j$-invariant is different from~$0$ and~$1728$.
The~key step in~\cite{P1} is to show~that
\begin{equation}
\label{g1n2}
n_{1,d}(\mu)=n_d^{(1)}(\mu),
\end{equation}
if $\mu$ is a tuple of $3d\!-\!1$ points in~$\PP$.
One of the main ingredients in proving Theorem~\ref{main_thm} is
Proposition~\ref{g1degen_prp}, which states 
that~\e_ref{g1n2} is valid for any tuple~$\mu$
that satisfies condition~\e_ref{main_thm_e}.
Note that the numbers listed in the above table are consistent
with~\e_ref{g1n2} and facts of classical algebraic geometry.
In particular, the image of every degree-$4$ map from a genus-one curve
to~$\P$ lies in a $\PPP$ and 
the image of every degree-$6$ map lies in a $\Bbb{P}^5$;
see~\cite[p116]{ACGH}.
Thus, the first three numbers in the table should be the same,
and the last two numbers should be the same.
The proof of Proposition~\ref{g1degen_prp}
extends the degeneration argument of~\cite{P1}
and builds up on modifications described in~\cite{Z1}.
We work with the moduli space $\ov{\frak M}_{1,N}(\P,d)$
of stable \hbox{degree-$d$} maps from genus-one $N$-pointed curves into~$\P$
and study what happens in the limit to the maps that pass through
the constraints~$\mu$ as the $j$-invariant of the domain tends infinity,
i.e.~the domain degenerates to a rational curve with two points identified.
\\

\noindent
Proposition~\ref{g1degen_prp} is not useful for determining
the numbers~$n_{1,d}(\mu)$ in $\P$ if $n\!\ge\!3$,
since the right-hand side of~\e_ref{g1n2} is unknown.
Computation of~$n_{1,d}(\mu)$ for all projective spaces
is the subject of~\cite{I},
where an entirely different approach is taken.
The main step in computing these numbers is showing that
$$2n_d^{(1)}(\mu)=\RT_{1,d}(\mu_1;\mu_2,\ldots,\mu_N)-CR_1(\mu),$$
where $CR_1(\mu)$ is the number of zeros of 
an explicit affine map between vector bundles 
over~$\bar{\cal V}_1'(\mu)$; see Proposition~\ref{g1_prp}.
The remaining step is to express this number of zeros topologically.
In general, if the linear part of an affine map~$\psi$ does not vanish,
it is easy to determine the signed cardinality of~$\psi^{-1}(0)$;
see Lemma~\ref{zeros_main}.
The approach of~\cite{I} is to replace
the linear part~$\al$ of the affine part under consideration
by a nonvanishing linear map over a space obtained 
from~$\bar{\cal V}_1'(\mu)$ by sequence of blowups
and then to express the resulting intersection number 
in terms of intersection numbers on the spaces~$\bar{\cal V}_k'(\mu)$.
The main problem with this approach is that the new linear map is
not described in~\cite{I} and it is not clear how to construct it in general.
In addition, the normal bundles of certain spaces needed 
for the second part of this approach are given incorrectly;
see Lemma~2.8 or equation~(2.27) in~\cite{I} for example.
Both of these statements can be corrected without affecting
the computability of the intersection numbers,
but presumably with a change in the final result.
If $n\!=\!2$, no blowup is needed.
If $n\!=\!3,4$, the zero set of~$\al$ is a complex manifold and
the ``derivative'' of~$\al$ in the normal direction 
along~$\al^{-1}(0)$ is nondegenerate.
In such cases, only one blowup is needed and 
a linear map with the required properties can be constructed fairly easily.
Furthermore, Lemma~2.8 of~\cite{I} requires no correction
in the $n\!=\!2,3,4$ cases, while
equation~(2.27) is never used.
If $n\!=\!2,3$, $CR_1(\mu)$ and $n_{1,d}(\mu)$ are then
expressed in terms of the numbers~$n_{d'}(\mu')$,
with $d'\!\le\! d$ and $\mu'$ related to~$\mu$.
Several numbers $n_{1,d}(\mu)$ for~$\Bbb{P}^4$ are given in~\cite{I}
as well.
However, 
no topological formula, like that of Theorem~\ref{main_thm},
is given for $CR_1(\mu)$ or $n_{1,d}(\mu)$ for $\P$ with $n\!\ge\!4$
and no number $n_{1,d}(\mu)$ is given for $\P$ \hbox{with $n\!\ge\!5$}.\\

\noindent
We obtain the expression of Theorem~\ref{main_thm}
for the number $CR_1(\mu)$ in Section~\ref{comp_sect};
see Proposition~\ref{cr_str}.
Our approach involves no blowups and requires relatively little understanding
of the global structure of the spaces~$\bar{\cal V}_k(\mu)$.
Instead we describe $CR_1(\mu)$ as the euler class of a bundle
minus the sum of contributions to the euler class from smooth, 
but usually noncompact, strata of the zero set of
the linear part~$\al_{1,0}$ of the affine map.
Computation of these contributions in good cases involves 
counting the zeros of affine maps again, but with the rank
of the target bundle reduced by one; see Subsection~\ref{top_sec}.
Of course, if we are to have any hope of computing these contributions, 
we need to understand the behavior of~$\al_{1,0}$
near the smooth strata of its zero~set.
Proposition~\ref{rat_str_prp} describes the behavior of~$\al_{1,0}$
and of related linear maps near the boundary strata 
of~$\bar{\cal V}_k(\mu)$.\\

\noindent
Theorem~\ref{main_thm} follows immediately from 
Propositions~\ref{cr_str} and~\ref{g1degen_prp}.
Their proofs are mutually independent.
Section~\ref{degen_sect} uses some of the notation defined 
in Subsection~\ref{notation}.
The topological tools of Subsection~\ref{top_sec},
the descriptive notation of Subsection~\ref{notation},
and the structure theorem of Subsection~\ref{str_sec}
are integral to the computations of Section~\ref{comp_sect}.\\

\noindent
In brief, we enumerate one-nodal rational curves from 
genus-one fixed-complex-structure invariants.
Can a similar approach be used with higher-genera enumerative invariants?
Let $\mu$ be an $N$-tuple of proper subvarieties
of~$\P$ in general position such that
$$\codim_{\Bbb{C}}\mu=d(n+1)-n.$$
Denote by $n_{2,d}(\mu)$ the number of genus-two degree-$d$
curves that pass through the constraints~$\mu$ and have a fixed 
generic complex structure on the normalization.
Let $n_d^{(3)}(\mu)$, $\tau_d(\mu)$, and $T_d(\mu)$ denote 
the number of rational two-component curves connected at three nodes, 
of rational curves with a triple point,
and of rational curves with a tacnode, respectively.
If $n\!=\!2$, we take $n_d^{(3)}(\mu)$ to be
 the number of two-component rational curves
with a choice of three nodes common to both components.
In all cases, the curves have degree-$d$ and pass through 
the constraints~$\mu$.
Completing the degeneration argument of~\cite{KQR},
it is shown in~\cite{Z1}~that
\begin{equation}
\label{g2n2}
n_{2,d}(\mu)=6\big(n_d^{(3)}(\mu)+\tau_d(\mu)+T_d(\mu)\big),
\end{equation}
if $\mu$ is a tuple of $3d\!-\!2$ points in~$\PP$.
The arguments of~\cite{KQR} and~\cite{Z1} should extend 
to show that equation~\e_ref{g2n2} is valid for 
arbitrary constraints~$\mu$ in all projective spaces.
On the other hand, $n_{2,d}(\mu)$ for~$\PPP$ is computed in~\cite{Z2}
and the method extends at least to~$\Bbb{P}^4$.
Thus, in those two cases, we should be able to express
the sum of the numbers $n_d^{(3)}(\mu)$, $\tau_d(\mu)$, and $T_d(\mu)$
in terms of intersection numbers of the spaces~$\bar{\cal V}_k(\mu)$.
The relation~\e_ref{g2n2} is obtained by considering a degeneration
to  a specific singular genus-two curve.
Perhaps, different relations can be obtained by considering
degeneration to other singular genus-two curves.
With enough different relations, we would be able to 
compute the numbers $n_d^{(3)}(\mu)$, $\tau_d(\mu)$, and $T_d(\mu)$
at least for $\PPP$ and~$\Bbb{P}^4$.\\

\noindent
The author thanks T.~Mrowka for many useful discussions
and E.~Ionel for comments on the original version of this paper.

\section{Background}

\subsection{Topology}
\label{top_sec}

\noindent
We begin by describing the topological tools used in the next section.
In particular, we review the notion of 
contribution to the euler class of a vector bundle from 
a (not necessarily closed) subset of the zero set of a section.
We also recall how one can enumerate the zeros of 
an affine map between vector bundles.
These concepts are closely intertwined.
Details can be found in Section~3 of~\cite{Z2}.\\

\noindent
Throughout this paper, all vector bundles 
are assumed to be complex and normed.
If $F\!\lra\!{\cal M}$ is a smooth vector bundle,
closed subset $Y$ of $F$ is {\it small}
if it contains no fiber of~$F$ and is preserved 
under scalar multiplication.
If ${\cal Z}$ is a compact oriented zero-dimensional manifold,
we denote the signed cardinality of~${\cal Z}$ by~$^{\pm}|{\cal Z}|$.
If $k$ is an integer, we write $[k]$ for the set of positive integers
not exceeding~$k$.

\begin{dfn}
\label{top_dfn1b}
Suppose $F,{\cal O}\!\lra\!{\cal M}$ are smooth vector bundles.\\
(1) If $F\!=\!\bigoplus\limits_{i=1}^{i=k}F_i$
and $\under{d}\!=\!(d_1,\ldots,d_k)$ is a $k$-tuple of positive integers,
bundle map $\al\!:F\!\lra\!{\cal O}$ 
is a \under{polynomial of degree~$\under{d}$} if
for each $i\!\in\![k]$ there~exists
$$p_i\!\in\!\Ga({\cal M};F_i^{*\otimes d_i}\!\otimes\!{\cal O})
\hbox{~~for~}i\!\in\![k]
\qquad\hbox{s.t.}\qquad
\al(\ups)=\sum_{i=1}^{i=k}p_i\big(\ups_i^{d_i}\big)
\quad\forall \ups=(\ups_i)_{i\in[k]}\in\bigoplus_{i=1}^{i=k}F_i.$$
(2) If $\al\!:F\!\lra\!{\cal O}$ is a polynomial, 
the \under{rank of~$\al$} is the number
$$\rk~\al\equiv\max\{\rk_b\al\!: b\!\in\!{\cal M}\},
\quad\hbox{where}\quad
\rk_b\al=\dim_{\Bbb{C}}\big(\hbox{Im}~\al_b\big).$$
Polynomial $\al\!:F\!\lra\!{\cal O}$ is \under{of constant rank}
if $\rk_b\al\!=\!\rk~\!\al$ for all $b\!\in\!{\cal M}$;
$\al$ is \under{nondegenerate}
if \hbox{$\rk_b\al=rk~\!F$} for all $b\!\in\!{\cal M}$.\\
(3) If $\Om$ is an open subset of $F$ and 
$\phi\!: \Om\!\lra\!{\cal O}$ is a smooth bundle map, 
bundle map $\al\!:F\!\lra\!{\cal O}$ is
a \under{dominant term of~$\phi$} 
if there exists $\ve\!\in\! C^0(F; \Bbb{R})$ such that
$$\big|\phi(\ups)-\al(\ups)\big|\le
\ve(\ups)|\al(\ups)|
\quad\forall\ups\!\in\!\Om
\quad\hbox{and}\quad
\lim_{\ups\lra0}\ve(\ups)=0.$$
Dominant term $\al\!:F\!\lra\!{\cal O}$ of $\phi$ is 
the \under{resolvent} of $\phi$ if $\al$ is a polynomial
of constant rank.
(4) $\phi\!: \Om\!\lra\!{\cal O}$ is \under{hollow}
if there exist dominant term $\al$ of $\phi$ 
and splittings $F\!=\! F^-\!\oplus F^+$
and ${\cal O}\!=\!{\cal O}^-\!\oplus {\cal O}^+$ 
such that $\al(F^+)\!\subset\!{\cal O}^+$,
$\al^-\!\equiv\!\pi^-\!\circ(\al|F^-)$ is a constant-rank polynomial,
where $\pi^-\!:{\cal O}\!\lra\!{\cal O}^-$ is the projection map,
and \hbox{$(\rk~\al^-\!+\!\frac{1}{2}\dim{\cal M})\!<\!\rk~{\cal O}^-$}.
\end{dfn}

\noindent
The base spaces we work with in the next two sections
are closely related to 
spaces of rational maps into $\P$ of total degree~$d$
that pass through the $N$ constraints $\mu_1,\ldots,\mu_N$.
From the algebraic geometry point of view,
spaces of rational maps are algebraic stacks, 
but with a fairly obscure local structure.
We view these spaces as {\it mostly smooth}, or {\it ms-}, manifolds:
compact oriented topological manifolds
stratified by smooth manifolds, such that the boundary strata
have (real) codimension at least two.
Subsection~\ref{str_sec} gives explicit descriptions
of neighborhoods of boundary strata and
of the behavior of certain bundle sections near such strata.
We call the main stratum ${\cal M}$ of ms-manifold $\bar{\cal M}$
the {\it smooth base} of~$\bar{\cal M}$.
Definition~3.7 in~\cite{Z2}
also introduces the natural notions of
{\it ms-maps} between ms-manifolds,
{\it ms-bundles} over ms-manifolds, and 
{\it ms-sections} of ms-bundles.

\begin{dfn}
\label{euler_dfn2}
Let  $\bar{\cal M}\!=\!{\cal M}_n\sqcup\bigsqcup_{i=0}^{n-2}{\cal M}_i
\!=\!{\cal M}\sqcup\bigsqcup_{i=0}^{n-2}{\cal M}_i$ 
be an ms-manifold of dimension~$n$.\\
(1) If ${\cal Z}\!\subset\!{\cal M}_i$ is a smooth oriented submanifold, 
a \under{normal-bundle model for ${\cal Z}$} is a tuple
$(F,Y,\vt)$, where\\
(1a) $F\!\lra\!{\cal Z}$ is a smooth vector bundle
and $Y$ is a small subset of~$F$;\\
(1b) for some $\de\!\in\!C^{\i}({\cal Z};\Bbb{R}^+)$,
$\vt\!:F_{\de}\!-\!(Y\!-\!{\cal Z})\lra\!\bar{\cal M}$ 
is a continuous map such that\\ 
(1b-i) $\vt\!:F_{\de}\!-\!(Y\!-\!{\cal Z})\!\lra\!\bar{\cal M}$ 
is a homeomorphism
onto an open neighborhood of ${\cal Z}$ in ${\cal M}\cup{\cal Z}$;\\
(1b-ii) $\vt|_{\cal Z}$ is the identity map, and
$\vt\!:F_{\de}\!-\!Y\!-\!{\cal Z}\!\lra\!{\cal M}$
is an orientation-preserving diffeomorphism on an open subset 
of~${\cal M}$.\\
(2) A \under{closure} of normal-bundle model $(F,Y,\vt)$ for ${\cal Z}$
is a tuple $(\bar{\cal Z},\tilde{F},\pi)$, where\\
(2a) $\bar{\cal Z}$ is an ms-manifold with 
smooth base ${\cal Z}$;\\
(2b) $\pi\!:\bar{\cal Z}\!\lra\!\bar{\cal M}$ is an ms-map
such that $\pi|_{\cal Z}$ is the identity;\\
(2c) $\tilde{F}\!\lra\!\bar{\cal Z}$ is an ms-bundle such that
$\tilde{F}|_{\cal Z}\!=\!F$.
\end{dfn}

\noindent
We use a normal-bundle model for ${\cal Z}$ to describe
the behavior of bundle sections over $\bar{\cal M}$ near~${\cal Z}$.
In particular, if $\al\!:E\!\lra\!{\cal O}$ is an ms-polynomial,
we call ${\cal Z}$ an {\it $\al$-regular} subset of $\bar{\cal M}$
if for some normal-bundle model $(F,Y,\vt)$ for ${\cal Z}$,
$\vt^*\al$ can be approximated, by a constant-rank polynomial
$p\!:F\oplus E\!\lra\!{\cal O}$;
see Definition~3.9 in~\cite{Z2}.
Polynomial $\al\!:E\!\lra\!{\cal O}$ is {\it regular}
if $\bar{\cal M}$ can be decomposed into finitely many
$\al$-regular subsets.
If $\rk~\!E\!+\!\frac{1}{2}\dim\bar{\cal M}\!=\!\rk~\!{\cal O}$,
for a generic $\bar{\nu}\!\in\!\Ga(\bar{\cal M};{\cal O})$,
the zero set of the polynomial map
$$\psi_{\al,\bar{\nu}}\!:E\lra{\cal O},\qquad 
\psi_{\al,\bar{\nu}}(\ups)=\nu_{\ups}+\al(\ups),$$
is a zero-dimensional oriented submanifold of~$E|{\cal M}$.
By Lemma~3.10 in~\cite{Z2},
if $\al$ is a regular polynomial,
$\psi_{\al,\bar{\nu}}^{-1}(0)$ is a finite set for a generic choice 
of~$\bar{\nu}$,
and $N(\al)\equiv^{\pm}\!\big|\psi_{\al,\bar{\nu}}^{-1}(0)\big|$
is independent of such a choice of~$\bar{\nu}$.\\

\noindent
As described below, counting the zeros of $\psi_{\al,\bar{\nu}}$
involves determining the contribution ${\cal C}_{\bar{\cal Z}}(s)$
to the euler class of a bundle~$V$
from a subset~$\bar{\cal Z}$ of the zero set of a section~$s$ of~$V$.
In the cases we encounter in Section~\ref{comp_sect}, 
$\bar{\cal Z}$~decomposes into disjoint, and usually non-compact, 
complex manifolds~${\cal Z}_i$ 
near which the behavior of~$s$ can be understood.
Then ${\cal C}_{\bar{\cal Z}}(s)\!=\!\sum{\cal C}_{{\cal Z}_i}(s)$,
where ${\cal C}_{{\cal Z}_i}(s)$ is
the \hbox{\it $s$-contribution of ${\cal Z}_i$ to $e(V)$}.
This is the signed number of elements of $\{s\!+\!\nu\}^{-1}(0)$
that lie very close to~${\cal Z}_i$,
where $\nu\!\in\!\Ga(\bar{\cal M};V)$ is a small generic perturbation of~$s$.
The manifolds ${\cal Z}_i$ we encounter fall in one of 
the two categories described below.

\begin{dfn}
\label{euler_dfn4}
Suppose $\bar{\cal M}$ is an ms-manifold of dimension~$2n$,
$V\!\lra\!\bar{\cal M}$ is an ms-bundle of rank~$n$,
$s\!\in\!\Ga(\bar{\cal M};V)$, and 
${\cal Z}\!\subset\!s^{-1}(0)$.\\
(1) ${\cal Z}$ is \under{$s$-hollow} 
if there exist a normal-bundle model $(F,Y,\vt)$ for ${\cal Z}$
and a bundle isomorphism 
\hbox{$\vt_V\!\!:\vt^*V\!\lra\!\pi_F^*V$},
covering the identity on $F_{\de}\!-\!(Y\!-\!{\cal Z})$, such~that\\
(1a) $\vt_V|_{F_{\de}-Y-{\cal Z}}$ is smooth and
$\vt_V|_{\cal Z}$ is the identity;\\
(1b) the map $\phi\equiv\vt_V\circ\vt^*s\!: 
          F_{\de}\!-\!(Y\!-\!{\cal Z})\!\lra\!V$ is hollow.\\
(2) ${\cal Z}$ is \under{$s$-regular} if there exist 
a normal-bundle model $(F,Y,\vt)$ for ${\cal Z}$
with closure~$(\bar{\cal Z},\tilde{F},\pi)$,
regular polynomial $\al\!:\tilde{F}\!\lra\!\pi^*V$, and 
a bundle isomorphism 
$\vt_V\!\!:\vt^*V\!\lra\!\pi_F^*V$
covering the identity on~$F_{\de}\!-\!(Y\!-\!{\cal Z})$, such~that\\
(2a) $\vt_V|_{F_{\de}-Y-{\cal Z}}$ is smooth and
$\vt_V|_{\cal Z}$ is the identity;\\
(2b)  $\al|_{\cal Z}$ is nondegenerate and is the resolvent for
$\phi\!\equiv\!\vt_V\!\circ\!\vt^*s\!: 
          F_{\de}\!-\!(Y\!-\!{\cal Z})\!\lra\!V$, 
and $Y$ is preserved under  scalar multiplication
in each of the components of $F$ for the splitting 
corresponding to $\al$ as in (1) of Definition~\ref{top_dfn1b}.
\end{dfn}

\begin{prp}
\label{euler_prp}
Let $V\!\lra\!\bar{\cal M}$ be an ms-bundle of rank~$n$
over an ms-manifold of dimension~$2n$.
Suppose ${\cal U}$ is an open subset of ${\cal M}$
and $s\!\in\!\Ga(\bar{\cal M};V)$ is such that
$s|_{\cal U}$ is transversal to the zero~set.\\
(1) If $s^{-1}(0)\cap{\cal U}$ is a finite set,
$^{\pm}|s^{-1}(0)\cap{\cal U}|
=\lan e(V),[\bar{\cal M}]\ran-{\cal C}_{\bar{\cal M}-{\cal U}}(s)$.\\
(2) If $\bar{\cal M}-{\cal U}=
  \bigsqcup\limits_{i=1}\limits^{i=k}\!{\cal Z}_i$,
where each ${\cal Z}_i$ is $s$-regular or $s$-hollow, then 
$s^{-1}(0)\cap{\cal U}$ is finite, and 
$$^{\pm}\big|s^{-1}(0)\cap{\cal U}\big|
=\big\lan e(V),[\bar{\cal M}]\big\ran-{\cal C}_{\bar{\cal M}-{\cal U}}(s)
=\big\lan e(V),[\bar{\cal M}]\big\ran-
\sum_{i=1}^{i=k}{\cal C}_{{\cal Z}_i}(s).$$
If ${\cal Z}_i$ is $s$-hollow, ${\cal C}_{{\cal Z}_i}(s)\!=\!0$.
If ${\cal Z}_i$ is $s$-regular 
and \hbox{$\al_i\!:\tilde{F}_i\!\lra\!V$} is the corresponding polynomial,
$${\cal C}_{{\cal Z}_i}(s)=
^{\pm}\!\big|\{\ups\!\in\!\tilde{F}_i\!:
 \bar{\nu}_{\ups}\!+\!\al_i(\ups)=0\}\big|
\equiv N(\al_i),$$
where $\bar{\nu}\!\in\!\Ga(\bar{\cal Z}_i;V)$ is a generic section.
Finally, if $\al_i\!\in\!
\Ga(\bar{\cal Z}_i;\tilde{F}_i^{*\otimes k}\otimes\pi^* V)$  
has constant rank  over~$\bar{\cal Z}_i$
and factors through a $\tilde{k}$-to-$1$ cover
$\rho_i\!:\tilde{F}_i\!\lra\!\tilde{F}_i^{\otimes k}$,
$${\cal C}_{{\cal Z}_i}(s)=\tilde{k}
\big\lan e\big(\pi^*V/\al_i(\tilde{F}_i)\big),[\bar{\cal Z}_i]\big\ran.$$
\end{prp}

\noindent
This is Corollary~3.13 in~\cite{Z2}.
Proposition~\ref{euler_prp} reduces the problem of 
computing ${\cal C}_{{\cal Z}_i}(s)$ for
an $s$-regular manifold ${\cal Z}_i$
to counting the zeros of a polynomial map between two vector bundles. 
The general setting for the latter problem is the following.
Suppose  $E,{\cal O}\!\lra\!\bar{\cal M}$ are ms-bundles such that
$\rk~E\!+\!\frac{1}{2}\dim\bar{\cal M}\!=\!\rk~{\cal O}$,
and $\al\!:E\!\lra\!{\cal O}$ is a regular polynomial.
Let $\bar{\nu}\!\in\!\Ga(\bar{\cal M};{\cal O})$
be such that the~map
$$\psi_{\al,\bar{\nu}}\!\equiv\!\bar{\nu}\!+\!\al\!: E\!\lra\!{\cal O}$$
is transversal to the zero set in ${\cal O}$ 
on~$E|{\cal M}$, and all its zeros are contained in~$E|{\cal M}$.
Then \hbox{$N(\al)\!\equiv^{\pm}\!|\psi_{\al,\bar{\nu}}^{-1}(0)|$}
depends only on~$\al$.
If the rank of $E$ is zero, then clearly
$$N(\al)=^{\pm}\!\big|\psi_{\al,\bar{\nu}}^{-1}(0)\big|
=\big\lan e({\cal O}),[\bar{\cal M}]\big\ran.$$
If the rank of $E$ is positive and $\bar{\nu}$ is generic,
it does not vanish and thus determines 
a trivial line subbundle $\Bbb{C}\bar{\nu}$ of~${\cal O}$.
Let ${\cal O}^{\perp}\!=\!{\cal O}/\Bbb{C}\bar{\nu}$
and denote by $\al^{\perp}$ the composition of $\al$ with
the quotient projection map.
If $E$ is a line bundle and $\al$ is~linear,
$$N(\al)=^{\pm}\!\big|\psi_{\al,\bar{\nu}}^{-1}(0)\big|
=\big\lan e(E^*\!\otimes\!{\cal O}^{\perp}),[\bar{\cal M}]\big\ran
-{\cal C}_{\al^{-1}(0)}(\al^{\perp}).$$
By Proposition~\ref{euler_prp}, computation of 
${\cal C}_{\al^{-1}(0)}(\al^{\perp})$ again involves counting
the zeros of polynomial maps, but with the rank of the new target bundle, 
i.e.~$E^*\!\otimes\!{\cal O}^{\perp}$, one less than the rank
of the original one, i.e.~${\cal O}$.
Subsection~3.3 in~\cite{Z2}
reduces the problem of determining $N(\al)$ in all other cases
to the case $E$ is a line bundle and $\al$ is linear.
Thus, at least in reasonably good cases, $N(\al)$ can be determined
after a finite number of steps.\\

\noindent
The next lemma summarizes the results of 
Subsection~3.3 in~\cite{Z2}
in the case the original map $\al\!:E\!\lra\!{\cal O}$ is linear.
This case suffices for our purposes.
We denote by 
$$\al'\in\Ga\big(\Bbb{P}E;\hbox{Hom}(\ga_E,\pi_E^*{\cal O})\big)$$
the section induced by~$\al$.
\hbox{Let $\la_E\!=\!c_1(\ga_E^*)$}.

\begin{lmm}
\label{zeros_main}
Suppose  $\bar{\cal M}$ is  an ms-manifold and
$E,{\cal O}\lra\bar{\cal M}$ are ms-bundles such that
$$\rk~E+\frac{1}{2}\dim\bar{\cal M}= \rk~{\cal O}.$$
If
$\al\!\in\!\Ga(\bar{\cal M};E^*\!\otimes\!{\cal O})$ and
$\bar{\nu}\!\in\!\Ga(\bar{\cal M};{\cal O})$
are such that $\al$ is regular,
$\bar{\nu}$ has no zeros, the map
$$\psi_{\al,\bar{\nu}}\!\equiv\!\bar{\nu}\!+\!\al\!: E\lra{\cal O}$$
is transversal to the zero set on~$E|{\cal M}$, and all its zeros
are contained in~$E|{\cal M}$, then
$\psi_{\al,\bar{\nu}}^{-1}(0)$ is a finite set,
$^{\pm}\!|\psi_{\al,\bar{\nu}}^{-1}(0)|$ depends only on $\al$, and
$$N(\al)\equiv ^{\pm}\!|\psi_{\al,\bar{\nu}}^{-1}(0)|
=\big\lan c({\cal O})c(E)^{-1}\!,[\bar{\cal M})]\big\ran
-{\cal C}_{{\al'}^{-1}(0)}({\al'}^{\perp}).$$
Furthermore, if $n\!=\!\rk~\!E$,
\begin{equation}\label{zeros_main_e}
\la_E^n+\sum_{k=1}^{k=n}c_k(E)\la_E^{n-k}=0\in H^{2n}(\Bbb{P}E)
\hbox{and}\quad
\big\lan \mu\la_E^{n-1},[\Bbb{P}E]\big\ran=
\big\lan\mu,[\bar{\cal M}]\big\ran~~
\forall\mu\!\in\!H^{2m-2n}(\bar{\cal M}).
\end{equation}
\end{lmm}

\subsection{Notation}
\label{notation}

\noindent
In this subsection, we describe the most important notation used 
in this paper. 
Some of the notation is only sketched;
see Section~2 in~\cite{Z3} for more details.\\

\noindent
If $I_1$ and $I_2$ are two sets, denote the disjoint
union of $I_1$ and $I_2$ by $I_1\!+\!I_2$.
We set 
$$\i=(0,0,-1)\in S^2\subset\Bbb{R}^3
\quad\hbox{and}\quad
e_{\i}\!=\!(1,0,0)\in T_{\i}S^2.$$
Let $q_N\!: \Bbb{C}\!\lra\! S^2\!\subset\!\Bbb{R}^3$ 
be the stereographic projection mapping the origin in $\Bbb{C}$ 
to the north pole.
We identify $\Bbb{C}$ with  $S^2\!-\!\{\i\}$ via the map~$q_N$.

\begin{dfn}
\label{rooted_tree}
A finite partially ordered set $I$ is a \under{linearly ordered set}
if for all \hbox{$i_1,i_2,h\!\in\! I$} such that $i_1,i_2\!<\!h$, 
either $i_1\!\le\! i_2$ \hbox{or $i_2\!\le\! i_1$.}\\
A linearly ordered set $I$ is a \under{rooted tree} if
$I$ has a unique minimal element, 
i.e.~there exists \hbox{$\hat{0}\!\in\! I$} such that $\hat{0}\!\le\! i$ 
for {all $i\!\in\!I$}.
\end{dfn}

\noindent
If $I$ is a linearly ordered set, let $\hat{I}$ be
the subset of the non-minimal elements of~$I$.
For every $h\!\in\!\hat{I}$,  denote by $\io_h\!\in\!I$
the largest element of $I$ which is smaller than~$h$.
Suppose $I\!=\!\bigsqcup\limits_{k\in K}\!I_k$
is the splitting of $I$ into rooted trees 
such that $k$ is the minimal element of~$I_k$.
If $\hat{1}\!\not\in\!I$, we define
the linearly ordered set $I\!+_k\!\hat{1}$ to
be the set $I\!+\!\hat{1}$ with all partial-order relations of~$I$
along with the relations
$k\!<\!\hat{1}$ and $\hat{1}\!<\!h$ \hbox{if $h\!\in\!\hat{I}_k$}.
\\

\noindent
If $S$ is a (possibly singular) complex curve and
$M$ is a finite set, 
a {\it $\P$-valued bubble map with $M$-marked 
points} is a tuple
$$b=\big(S,M,I;x,(j,y),u\big),$$
where $I$ is a linearly ordered set, and
$$x\!:\hat{I}\!\lra\!S\cup S^2,~~j\!:M\!\lra\! I,~~ 
y\!:M\!\lra\!S\cup S^2,\hbox{~~and~~} 
u\!:I\!\lra\!C^{\i}(S;\P)\cup C^{\i}(S^2;\P)$$
are maps such that
$$x_h\in
\begin{cases}
S^2\!-\!\{\i\},&\hbox{if~}\io_h\!\in\!\hat{I};\\
S,&\hbox{if~}\io_h\!\not\in\!\hat{I};
\end{cases}\qquad
y_l\in
\begin{cases}
S^2\!-\!\{\i\},&\hbox{if~}j_l\!\in\!\hat{I};\\
S,&\hbox{if~}j_l\!\not\in\!\hat{I};\qquad
\end{cases}
u_i\in
\begin{cases}
C^{\i}(S^2;\P),&\hbox{if~}i\!\in\!\hat{I};\\
C^{\i}(S;\P),&\hbox{if~}i\!\not\in\!\hat{I};
\end{cases}$$  
and  $u_h(\i)\!=\!u_{\io_h}(x_h)$ for all 
$h\!\in\!\hat{I}$.
We associate such a tuple with Riemann surface
$$\Si_b=
\Big(\bigsqcup_{i\in I}\Si_{b,i}\Big)\Big/\!\sim,
\hbox{~~where}\qquad 
\Si_{b,i}=
\begin{cases}
\{i\}\!\times\! S^2,&\hbox{if~}i\!\in\!\hat{I};\\
\{i\}\!\times\! S,&\hbox{if~}i\!\not\in\!\hat{I},
\end{cases}
\quad\hbox{and}\quad
(h,\i)\sim (\io_h,x_h)
~~\forall h\!\in\!\hat{I},$$
with marked points $(j_l,y_l)\!\in\!\Si_{b,j_l}$,
and continuous map $u_b\!:\Si_b\!\lra\!\P$,
given by $u_b|\Si_{b,i}\!=\!u_i$ for \hbox{all $i\!\in\! I$}.
We require that all the singular points of $\Si_b$
and all the marked points be distinct.
Furthermore, if $S\!=\!S^2$, all these points are to be different
from the special marked point $(\hat{0},\i)\!\in\!\Si_{b,\hat{0}}$.
In addition, if $\Si_{b,i}\!=\!S^2$ and 
$u_{i*}[S^2]\!=\!0\!\in\! H_2(\P;\Bbb{Z})$,
then $\Si_{b,i}$ must contain at least two singular and/or
marked points of~$\Si_b$ other than~$(i,\i)$.
If $S\!\neq\!S^2$, but $S$ is unstable, 
$u_i$ must satisfy a similar stability condition whenever 
$\Si_{b,i}\!=\!S$.
In particular, if $S$ is a torus or a circle of spheres and 
the restriction of $u_i$ to a component~$S_h$ of~$S$ is 
homologically zero,
$S_h$~contains at least one marked point of~$\Si_b$.
Two bubble maps $b$ and $b'$ are {\it equivalent} if there
exists a homeomorphism $\phi\!:\Si_b\!\lra\!\Si_{b'}$
such that $u_b\!=\!u_{b'}\!\circ\phi$,
$\phi(j_l,y_l)\!=\!(j_l',y_l')$ for \hbox{all $l\!\in\! M$} and
$\phi|_{\Si_{b,i}}$ is holomorphic for \hbox{all $i\!\in\!I$}.\\

\noindent
The general structure of bubble maps is described
by tuples ${\cal T}\!=\!(S,M,I;j,\under{d})$,
with \hbox{$d_i\!\in\!\Bbb{Z}$} 
specifying the degree of the map $u_b$ on~$\Si_{b,i}$.
We call such tuples {\it bubble types}.
Bubble type ${\cal T}$ is {\it simple} if  $I$ is a rooted tree;
${\cal T}$ is {\it basic} if $\hat{I}\!=\!\eset$ and $d_i\!\neq\!0$
for all $i\!\in\!I$;
${\cal T}$ is {\it semiprimitive} if 
$\io_h\!\not\in\!\hat{I}$, \hbox{$d_{\io_h}\!=\!0$}, and $d_h\!\neq\!0$,
for all $h\!\in\!\hat{I}$.
The above equivalence relation on the set of bubble maps
induces an equivalence relation on the set of bubble types.
For each $h,i\!\in\! I$, let
\begin{gather*}
D_i{\cal T}=\{h\!\in\!\hat{I}\!:i\!<\!h\},\quad
\bar{D}_i{\cal T}=D_i{\cal T}\cup\{i\},\quad
H_i{\cal T}=\{h\!\in\!\hat{I}\!:\io_h\!=\!i\},\quad
M_i{\cal T}=\{l\!\in\! M\!:j_l\!=\!i\},\\
\chi_{\cal T}h=\begin{cases}
0,&\hbox{if~} \forall i\!\in\! 
             I\hbox{~s.t.~}h\!\in\!\bar{D}_i{\cal T}, d_i\!=\!0;\\
1,&\hbox{if~} d_h\!\neq\! 0,\hbox{~but~}
\forall i\!\in\! I\hbox{~s.t.~}h\!\in\! D_i{\cal T}, d_i\!=\!0;\\
2,&\hbox{otherwise};
\end{cases}\qquad
\chi({\cal T})=\big\{h\!\in\!I\!:\chi_{\cal T}h\!=\!1\big\}.
\end{gather*}
Denote by ${\cal H}_{\cal T}$ the space of all holomorphic
bubble maps with structure~${\cal T}$.\\

\noindent
The automorphism group of every bubble type ${\cal T}$ we encounter
in the next two sections is trivial.
Thus, every bubble type discussed below is presumed to be 
automorphism-free.\\

\noindent
If $S$ is a circle of spheres, we denote by ${\cal M}_{\cal T}$ 
the set of equivalence classes of bubble maps in~${\cal H}_{\cal T}$.
For each bubble type ${\cal T}=(S^2,M,I;j,\under{d})$, let
$${\cal U}_{\cal T}=\big\{[b]\!:
b\!=\!\big(S^2,M,I;x,(j,y),u\big)\!\in\!{\cal H}_{\cal T},~
u_{i_1}(\i)=u_{i_2}(\i)~\forall i_1,i_2\!\in\! I\!-\!\hat{I}
\big\}.$$
Then there exists ${\cal B}_{\cal T}\!\subset\!{\cal H}_{\cal T}$
such that 
${\cal U}_{\cal T}$ is the quotient of a subset ${\cal B}_{\cal T}$
of ${\cal H}_{\cal T}$ by a
\hbox{$\tilde{G}_{\cal T}\!\equiv\!(S^1)^I$}-action.
Denote by ${\cal U}_{\cal T}^{(0)}$ the quotient of ${\cal B}_{\cal T}$
by 
\hbox{$G_{\cal T}\!\equiv\!(S^1)^{\hat{I}}\!\subset\!\tilde{G}_{\cal T}$}.
Then ${\cal U}_{\cal T}$ is the quotient of ${\cal U}_{\cal T}^{(0)}$ 
by the residual 
$G_{\cal T}^*\!\equiv\!(S^1)^{I-\hat{I}}\!\subset\!\tilde{G}_{\cal T}$
action.
Corresponding to these quotients, we obtain line orbi-bundles
$\big\{L_i{\cal T}\!\!\lra\!{\cal U}_{\cal T}\!\!:i\!\in\! I\big\}$.
Let 
$${\cal FT}=\bigoplus_{h\in\hat{I}}{\cal F}_h{\cal T}\lra{\cal U}_{\cal T},
\hbox{~~where}\quad
{\cal F}_h{\cal T}=L_h{\cal T}\otimes L_{\io_h}^*{\cal T}.$$
Denote by ${\cal F}^{\eset}{\cal T}$ the open subset of
${\cal FT}$ consisting of vectors with all components nonzero.\\

\noindent
Gromov topology on the space of equivalence classes of bubble maps
induces a partial ordering on the set of bubble types
and their equivalence classes such that the spaces
$$\bar{\cal U}_{\cal T}^{(0)}=
\bigcup_{{\cal T}'\le{\cal T}}{\cal U}_{{\cal T}'}^{(0)}
\qquad\hbox{and}\qquad
\bar{\cal U}_{\cal T}=
\bigcup_{{\cal T}'\le{\cal T}}{\cal U}_{{\cal T}'}$$
are compact and Hausdorff.
The $G_{\cal T}^*$-action on ${\cal U}_{\cal T}^{(0)}$ extends
to an action on $\bar{\cal U}_{\cal T}^{(0)}$,
and thus the line orbi-bundles $L_i{\cal T}\!\lra\!{\cal U}_{\cal T}$
with $i\!\in\!I\!-\!\hat{I}$ extend over~$\bar{\cal U}_{\cal T}$.
These bundles can be identified with 
the universal tangent line bundles for appropriate sections
of the universal bundle over~$\bar{\cal U}_{\cal T}$.
The evaluation maps
$$\ev_l\!:{\cal H}_{\cal T}\lra\P,\quad
\ev_l\big((S,M,I;x,(j,y),u)\big)=
u_{j_l}(y_l),$$
descend to all the quotients and induce continuous
maps on $\bar{\cal U}_{\cal T}$ and~$\bar{\cal U}_{\cal T}^{(0)}$.
If $\mu\!=\!\mu_M$ is an $M$-tuple of subvarieties of~$\P$,
let
$${\cal M}_{\cal T}(\mu)=
\big\{b\!\in\!{\cal M}_{\cal T}\!:\ev_l(b)\!\in\!\mu_l~\forall 
l\!\in\! M\big\}$$
and define spaces ${\cal U}_{\cal T}(\mu)$,
$\bar{\cal U}_{\cal T}(\mu)$, etc.~in a similar way.
If $S=S^2$, we define another evaluation~map,
$$\ev\!: 
{\cal B}_{\cal T}\lra\P\quad\hbox{by}\quad
\ev\big((S^2,M,I;x,(j,y),u)\big)=u_{\hat{0}}(\i),$$
where $\hat{0}$ is any minimal element of $I$.
This map descends to ${\cal U}_{\cal T}^{(0)}$ and~${\cal U}_{\cal T}$.
If $\mu\!=\!\mu_{\tilde{M}}$ is an $\tilde{M}$-tuple of constraints, let
$${\cal U}_{\cal T}(\mu)=
\big\{b\!\in\!{\cal U}_{\cal T}\!:
\ev_l(b)\!\in\!\mu_l~\forall l\!\in\!M\cap\tilde{M},~
\ev(b)\!\in\!\mu_l~\forall l\!\in\!M\!-\!\tilde{M}\big\}$$
and define ${\cal U}_{\cal T}^{(0)}(\mu)$, etc.~similarly.\\

\noindent
Suppose ${\cal T}\!=\!(S^2,M,I;j,\under{d})$ is a bubble type, 
$k\!\in\! I\!-\!\hat{I}$, and $M_0$ is nonempty subset of~$M_k{\cal T}$.
Let
$${\cal T}/M_0=\big(S^2,I,M\!-\!M_0;
j|(M\!-\!M_0),\under{d}\big).$$
Define ${\cal T}(M_0)\equiv(S^2,M,I+_k\hat{1};j',\under{d}')$  by
$$j'_l=\begin{cases}
k,&\hbox{if~}l\!\in\! M_0;\\
\hat{1},&\hbox{if~}l\!\in\! M_k{\cal T}\!-\!M_0;\\
j_l,&\hbox{otherwise};
\end{cases}\qquad
d'_i=\begin{cases}
0,&\hbox{if~}i\!=\!k;\\
d_k,&\hbox{if~}i\!=\!\hat{1};\\
d_i,&\hbox{otherwise}.
\end{cases}$$
The tuples ${\cal T}/M_0$ and ${\cal T}(M_0)$ are bubble types as long as 
$d_k\!\neq\!0$ or $M_0\!\neq\! M_{\hat{0}}{\cal T}$.
Then,
\begin{equation}
\label{cart_split2}
\bar{\cal U}_{{\cal T}(M_0)}(\mu)=
 \bar{\cal M}_{0,\{\hat{1}\}+M_0}\times \bar{\cal U}_{{\cal T}/M_0}(\mu),
\end{equation}
where $\bar{\cal M}_{0,\{\hat{1}\}+M_0}$ denotes the
Deligne-Mumford moduli space of rational curves with 
$(\{\hat{0},\hat{1}\}+M_0)$-marked points. 
If ${\cal T}$ is a basic bubble type, let
\begin{equation}\label{normal_bundle1}
c_1({\cal L}_k^*{\cal T})\equiv c_1(L_k^*{\cal T})-
\!\sum_{\eset\neq M_0\subset M_k{\cal T}}\!\!\!\!
PD_{\bar{\cal U}_{\cal T}(\mu)}\big[\bar{\cal U}_{{\cal T}(M_0)}(\mu)\big]
\in H^2\big(\bar{\cal U}_{\cal T}(\mu)\big).
\end{equation}
This cohomology class is well-defined;
see~Subsection~5.2 in~\cite{Z2}.\\

\noindent
We are now ready to explain the claim of Theorem~\ref{main_thm}.
Let $n$, $d$, $N$ and $\mu$ be as in the statement of the theorem.
If $k\!\ge\!1$ and $m\!\ge\!1$,
denote by $\bar{\cal V}_{k,m}(\mu)$ 
the disjoint union of the spaces $\bar{\cal U}_{\cal T}(\mu)$
taken over equivalence classes of basic bubble types
\hbox{${\cal T}\!=\!(S^2,[N]\!-\!M_0,I;j,\under{d})$} 
with $|M_0|\!=\!m$, $|I|\!=\!k$, and \hbox{$\sum d_k\!=\!d$}.
Let $\bar{\cal V}_k(\mu)\!=\!\bar{\cal V}_{k,0}(\mu)$.
We define the spaces ${\cal V}_{k,m}(\mu)$ similarly.
Let 
$$\big\{c_1({\cal L}_i^*)\!:i\!\in\![k]\big\},
\big\{c_1(L_i^*)\!:i\!\in\![k]\big\}\subset H^2
\big(\bar{\cal V}_{k,m}(\mu);\Bbb{Z}\big)$$ 
be given~by
$$\big\{c_1({\cal L}_i^*)\big|\bar{\cal U}_{\cal T}(\mu)\!:i\!\in\![k]\big\}
=\big\{c_1({\cal L}_i^*{\cal T})\!: i\!\in\!I\big\},
\quad
\big\{c_1(L_i^*)\big|\bar{\cal U}_{\cal T}(\mu)\!:i\!\in\![k]\big\}
=\big\{c_1(L_i^*{\cal T}):i\!\in\!I\big\},$$
whenever ${\cal T}$ is as above.
We denote by $\eta_l,\tilde{\eta}_l\!\in\!
H^{2l}\big(\bar{\cal V}_{k,m}(\mu);\Bbb{Z}\big)$ 
the sum of all degree-$l$ monomials in 
$\big\{c_1({\cal L}_i^*)\!:i\!\in\![k]\big\}$ and in
$\big\{c_1(L_i^*)\!:i\!\in\![k]\big\}$, respectively.
For example, 
$$\eta_2=c_1^2({\cal L}_1^*)+c_1^2({\cal L}_2^*)+
c_1({\cal L}_1^*)c_1({\cal L}_2^*)\in 
H^4\big(\bar{\cal V}_{k,m}(\mu);\Bbb{Z}\big).$$
Finally, let $a\!=\!\ev^*c_1(\ga_{\P}^*)\!\in\!
H^2\big(\bar{\cal V}_{k,m}(\mu);\Bbb{Z}\big)$, where 
$\ga_{\P}\!\lra\!\P$ denotes the tautological line bundle.\\

\noindent
We next describe a generalization of the splitting~\e_ref{cart_split2}
which is used in computations in Section~\ref{comp_sect}.
If \hbox{${\cal T}\!=\!(S^2,I,[N]\!-\!M_0;j,\under{d})$} is a bubble type,
let 
$$\bar{\cal T}\!=\!\big(S^2,\bar{I},[N]\!-\!\bar{M}_0;
j|([N]\!-\!\bar{M}_0),\under{d}|\bar{I}\big),
~~~\hbox{where}~~~
\bar{I}=I\!-\!\big\{i\!\in\!I\!-\!\hat{I}\!: d_i\!=\!0\big\},~~
\bar{M}_0=M_0\cup\!\bigcup_{i\in I-\bar{I}}\!\!\!M_i{\cal T}.$$
Note that if ${\cal T}$ is semiprimitive, $\bar{\cal T}$ is basic.
Furthermore,
\begin{gather}\label{cart_split1a}
{\cal U}_{\cal T}(\mu)=
\prod_{i\in I-\bar{I}}\!{\cal M}_{0,H_i{\cal T}+M_i{\cal T}}
\times {\cal U}_{\bar{\cal T}}(\mu),\\
\label{cart_split1b}
\bar{\cal U}_{\cal T}(\mu)=
\prod_{i\in I-\bar{I}}\!\bar{\cal M}_{0,H_i{\cal T}+M_i{\cal T}}\times
\bar{\cal U}_{\bar{\cal T}}(\mu),
\end{gather}
where ${\cal M}_{0,H_i{\cal T}+M_i{\cal T}}$ denotes the main
stratum of~$\bar{\cal M}_{0,H_i{\cal T}+M_i{\cal T}}$.
If $i\!\in\!I\!-\!\bar{I}$, by definition, the bundle 
\hbox{$L_i{\cal T}\!\lra\!\bar{\cal U}_{\cal T}(\mu)$}
is the pullback by the projection map of the bundle
$$L_{\hat{0}}{\cal T}_i^{(0)}\lra
\bar{\cal M}_{0,H_i{\cal T}+M_i{\cal T}}=
\bar{\cal U}_{{\cal T}_i^{(0)}},
\quad\hbox{where}\quad
{\cal T}_i^{(0)}=
\big(S^2,H_i{\cal T}\!+\!M_i{\cal T},\{\hat{0}\};\hat{0},0\big).$$
We call the latter bundle the {\it tautological line bundle}
over~$\bar{\cal M}_{0,H_i{\cal T}+M_i{\cal T}}$.
This is the universal tangent line at the marked point
$\hat{0}\!\in\!\bar{\cal M}_{0,H_i{\cal T}+M_i{\cal T}}$.\\

\noindent
Finally, if $X$ is any space, $F\!\lra\! X$ is a normed vector bundle,
and $\de\!: X\!\lra\!\Bbb{R}$ is any function,~let
$$F_{\de}=\big\{(b,v)\!\in\! F\!: |v|_b<\de(b)\big\}.$$
Similarly, if $\Om$ is a subset of $F$, let 
$\Om_{\de}=F_{\de}\cap\Om$.
If $\ups\!=\!(b,v)\!\in\! F$, denote by $b_{\ups}$ the image of $\ups$
under the bundle projection map, i.e.~$b$ in this case.

\subsection{A Structural Description}
\label{str_sec}

\noindent
We now describe the structure of the spaces~$\bar{\cal V}_{k,m}(\mu)$
and the behavior of certain bundle sections over~$\bar{\cal V}_{k,m}(\mu)$  
near the boundary strata.\\

\noindent
If $b\!=\!\big(S^2,M,I;x,(j,y),u\big)\!\in\!{\cal B}_{\cal T}$
 and $k\!\in\! I$, let
$${\cal D}_{{\cal T},k}b=du_k\big|_{\i}e_{\i}.$$
If ${\cal T}^*$ is a basic bubble type, 
the maps ${\cal D}_{{\cal T},k}$ with 
${\cal T}\!<\!{\cal T}^*$ and $k\!\in\!I\!-\!\hat{I}$
induce a continuous section of $\ev^*T\P$ over 
$\bar{\cal U}_{{\cal T}^*}^{(0)}$ and a continuous section of the bundle 
$L_k^*{\cal T}^*\!\otimes\!\ev^*T\P$ over 
$\bar{\cal U}_{{\cal T}^*}$, described~by 
$${\cal D}_{{\cal T}^*,k}[b,c_k]
=c_k{\cal D}_{{\cal T},k}b,
\quad\hbox{if}~~b\!\in\!{\cal U}_{\cal T}^{(0)},~c_k\!\in\!\Bbb{C}.$$

\begin{prp}
\label{rat_str_prp}
Suppose $p\!>\!2$, $n\!\ge\!2$, $d\!\ge\!1$, $N\!\ge\!1$,  
$\mu\!=\!(\mu_1,\ldots,\mu_N)$ is an $N$-tuple of proper subvarieties of~$\P$ 
in general position, such that
$$\codim_{\Bbb{C}}\mu\equiv\sum_{l=1}^{l=N}\codim_{\Bbb{C}}\mu_l-N
=d(n+1)-1,$$
and $M_0$ is a subset of~$[N]$.
If ${\cal T}^*\!\!=\!(S^2,[N]\!-\!M_0,I^*;j^*,\under{d}^*)$ is a 
basic bubble type such that \hbox{$\sum\! d_i^*\!=\!d$},
the space $\bar{\cal U}_{{\cal T}^*}(\mu)$ 
is an  ms-manifold of (real) dimension 
$2\big(n\!+\!1\!-\!2|I^*|\!-\!|M_0|\big)$
and $L_k{\cal T}^*$ for $k\!\in\!I^*$ and $\ev^*T\P$
are ms-bundles over~$\bar{\cal U}_{{\cal T}^*}(\mu)$.
If ${\cal T}\!=\!(S^2,[N]\!-\!M_0,I;j,\under{d})\!<\!{\cal T}^*$, 
there exist 
\hbox{$\de,C\!\in\! C^{\i}\big({\cal U}_{\cal T}(\mu);\Bbb{R}^+\big)$}
and a homeomorphism
$$\ga_{\cal T}^{\mu}\!: 
     {\cal FT}_{\de}\lra \bar{\cal U}_{{\cal T}^*}(\mu),$$
onto an open neighborhood of ${\cal U}_{\cal T}(\mu)$
in $\bar{\cal U}_{{\cal T}^*}(\mu)$ such that
$\ga_{\cal T}^{\mu}|{\cal U}_{\cal T}(\mu)$ is the identity
and $\ga_{\cal T}^{\mu}|{\cal F}^{\eset}{\cal T}_{\de}$
is an orientation-preserving diffeomorphism onto an open subset 
of~${\cal U}_{\cal T}(\mu)$.
Furthermore, with appropriate identifications,
\begin{gather*}
\Big|{\cal D}_{{\cal T}^*,k}\ga_{\cal T}^{\mu}(\ups)
-\al_{{\cal T},k}\big(\rho_{\cal T}(\ups)\big)\Big|
\le C(b_{\ups})|\ups|^{\frac{1}{p}} \big|\rho_{\cal T}(\ups)\big|
\quad\forall\ups\!\in\!{\cal FT}_{\de},
\qquad\hbox{where}\\
\rho_{\cal T}(\ups)=
\big((\tilde{\ups}_h)_{h\in\chi({\cal T})}\big)
\in\tilde{\cal F}{\cal T}\equiv\!
\bigoplus_{h\in\chi({\cal T})}\!\!
L_h{\cal T}\!\otimes\! L_{\tilde{\io}_h}^*{\cal T};~~
\tilde{\ups}_h=\!\bigotimes_{i\in\hat{I},h\in\bar{D}_i{\cal T}}
\!\!\!\!\ups_i;~~
\tilde{\io}_h\!\in\!I\!-\!\hat{I},~ 
h\!\in\!\bar{D}_{\tilde{\io}_h}{\cal T}\\
\al_{{\cal T},k}\big((\tilde{\ups}_h)_{h\in\chi({\cal T})}\big)
=\sum_{h\in I_k\cap\chi({\cal T})}\!\!\!\!\!
{\cal D}_{{\cal T},h}\tilde{\ups}_h,
\end{gather*}
and $I_k\!\subset\!I$ is the rooted tree containing~$k$.
\end{prp}

\noindent
This is a special case of Theorem~2.8 in~\cite{Z2};
see also the remark following the theorem.
The dimension of~$\bar{\cal U}_{{\cal T}^*}(\mu)$ is obtained 
as~follows:
\begin{equation*}\begin{split}
\frac{1}{2}\dim\bar{\cal U}_{{\cal T}^*}(\mu)=
\dim_{\Bbb{C}}{\cal U}_{{\cal T}^*}(\mu)
&=\sum_{i\in I^*}\big(d_i(n\!+\!1)+n-2\big)-(|I^*|-1)n-
\big(\codim_{\Bbb{C}}\mu+|M_0|\big)\\
&=n+1-2|I^*|-|M_0|.
\end{split}\end{equation*}
The analytic estimate on~${\cal D}_{{\cal T}^*,k}$
is crucial for implementation of the topological tools 
of Subsection~\ref{top_sec} in Subsection~\ref{comp_top}.
If ${\cal T}$ is semiprimitive,
the bundle ${\cal FT}\!=\!\tilde{\cal F}{\cal T}$
and the section $\al_{\cal T}\!=\al_{\cal T}\!\circ\!\rho_{\cal T}$
extend over $\bar{\cal U}_{\cal T}(\mu)$
via the decomposition~\e_ref{cart_split1b}.
In terms of the notions of Subsection~\ref{top_sec},
$({\cal FT},{\cal FT}\!-\!F^{\eset}{\cal T},\ga_{\cal T}^{\mu})$
is a normal-bundle model for 
${\cal U}_{\cal T}(\mu)\subset\bar{\cal U}_{{\cal T}^*}(\mu)$.
This normal-bundle model admits a closure if ${\cal T}$ is semiprimitive.
Note that ${\cal FT}$ is not usually the normal bundle 
of $\bar{\cal U}_{\cal T}(\mu)$ in $\bar{\cal U}_{{\cal T}^*}(\mu)$
if both spaces are viewed as algebraic stacks; see~\cite{P2}.
Proposition~\ref{rat_str_prp} implies only that  the restrictions 
to ${\cal U}_{\cal T}(\mu)$
of~${\cal FT}$ and of the normal bundle of 
$\bar{\cal U}_{\cal T}(\mu)$ in $\bar{\cal U}_{{\cal T}^*}(\mu)$
are isomorphic as topological vector bundles.\\

\noindent
For any $k,m\!\in\!\Bbb{Z}$, we define 
bundle $E_{k,m}\!\lra\!\bar{\cal V}_{k,m}(\mu)$
and homomorphism \hbox{$\al_{k,m}\!:E_{k,m}\!\lra\!\ev^*T\P$} 
over $\bar{\cal V}_{k,m}(\mu)$ by
$$E_{k,m}|\bar{\cal U}_{{\cal T}^*}(\mu)=
\bigoplus_{i\in I^*}L_i{\cal T}^*,\quad
\al_{k,m}\big((\ups_i)_{i\in I^*}\big)
=\sum_{i\in I^*}{\cal D}_{{\cal T}^*,i}\ups_i,$$
whenever ${\cal T}^*\!=\!(S^2,[N]\!-\!M_0,I^*;j^*,\under{d}^*)$ is a 
basic bubble type such that~$\sum\! d_i^*\!=\!d$, 
$|I^*|\!=\!k$, and $|M_0|\!=\!m$.
The following lemma will be used in Section~\ref{comp_sect}.

\begin{lmm}
\label{str_lmm}
Suppose $n\!\ge\!2$, $d\!\ge\!1$, $N\!\ge\!1$,  and
$\mu=(\mu_1,\ldots,\mu_N)$ is an $N$-tuple of proper subvarieties of~$\P$ 
in general position such that
$\codim_{\Bbb{C}}\mu\!=\!d(n\!+\!1)\!-\!1$.
If \hbox{${\cal T}\!=\!(S^2,[N]\!-\!M_0,I;j,\under{d})$} is a bubble type
such that ${\cal U}_{\cal T}(\mu)\!\subset\!\bar{\cal V}_{k,m}(\mu)$,
the restriction of~$\al_{k,m}$ to the~subbundle
$$E{\cal T}^{\perp}
\equiv\bigoplus_{i\in\chi({\cal T})-\hat{I}}\!\!\!\!L_i{\cal T}
\subset E_{k,m}$$
is nondegenerate over~${\cal U}_{\cal T}(\mu)$.
\end{lmm}

\noindent
{\it Proof:} The linear map~$\al_{k,m}$ has full rank on $E{\cal T}^{\perp}$
over ${\cal U}_{\cal T}(\mu)$ if and only if the section
$$\big\{\al_{k,m}|E{\cal T}^{\perp}\big\}'
\in\Ga\big(\Bbb{P}E{\cal T}^{\perp}|{\cal U}_{\cal T}(\mu);
\ga_{E{\cal T}^{\perp}}^*\!\otimes\!\ev^*T\P\big)$$ 
has no zeros. 
Note that
$$\dim_{\Bbb{C}} \Bbb{P}E{\cal T}^{\perp}|{\cal U}_{\cal T}(\mu)
\le \dim_{\Bbb{C}}{\cal V}_k(\mu)+(k-1) =n-k<n.$$
Thus, it is enough to show that $\big\{\al_{k,m}|E{\cal T}^{\perp}\big\}'$ 
is transversal to
the zero set in $\Bbb{P}E{\cal T}^{\perp}|{\cal U}_{\cal T}(\mu)$
if the constraints~$\mu$ are in general position.
This last fact is immediate from Lemma~\ref{reg_crl2}.

\begin{lmm}
\label{reg_crl2}
If $u\!:S^2\!\lra\!\P$ is a holomorphic map of positive degree
and $e_{\i}\!\in\!T_{\i}S^2$ is a nonzero vector,
the linear maps
\begin{alignat*}{2}
&H^0_{\bar{\partial}}(S^2;u^*T\P)\lra T_{u(\i)}\P,&\qquad&
\xi\lra\xi(\i),\\
&\big\{\xi\!\in\!H^0_{\bar{\partial}}(S^2;u^*T\P)\!:\xi(\i)\!=\!0\big\}
\lra T_{u(\i)}\P,
&\qquad &\xi\lra \nabla_{e_{\i}}\xi,
\end{alignat*}
are onto.
\end{lmm}

\noindent
This lemma is well-known; see
Corollary~6.3 in~\cite{Z2} for example.

\section{Computations}
\label{comp_sect}

\subsection{Topology}
\label{comp_top}

\noindent
In this section, we prove

\begin{prp}
\label{cr_str}
Suppose $n\!\ge\!2$, $d\!\ge\!1$, 
and $\mu\!=\!(\mu_1,\ldots,\mu_N)$ is 
an $N$-tuple of proper subvarieties of~$\P$ in general position
such~that
$$\codim_{\Bbb{C}}\mu\equiv
\sum_{l=1}^{l=N}\codim_{\Bbb{C}}\mu_l-N=d(n+1)-1.$$
Then the number of degree-$d$ genus-one curves that have 
a fixed generic complex structure on the normalization and 
pass through the constraints~$\mu$ is given~by
\begin{gather*}
n_{1,d}(\mu)=\frac{1}{2}\big(
\RT_{1,d}(\mu_1;\mu_2,\ldots,\mu_N)-CR_1(\mu)\big),
\qquad\hbox{where}\\
CR_1(\mu)=\sum_{k=1}^{2k\le n+1}\!\! (-1)^{k-1}(k\!-\!1)!
\!\sum_{l=0}^{n+1-2k}\!\! \binom{n\!+\!1}{l}
\big\lan a^l\eta_{n+1-2k-l},\big[\bar{\cal V}_k(\mu)\big]\big\ran.
\end{gather*}\\
\end{prp}

\noindent
We use the topological tools of Subsection~\ref{top_sec}
and the analytic estimate of Proposition~\ref{rat_str_prp}
to deduce Proposition~\ref{cr_str} from Proposition~\ref{g1_prp}.
The main step is Lemma~\ref{comp_top1}; 
the rest of this section is essentially combinatorics.

\begin{prp}
\label{g1_prp}
Suppose $n\!\ge\!2$, $d\!\ge\!1$, 
and $\mu\!=\!(\mu_1,\ldots,\mu_N)$ is 
an $N$-tuple of proper subvarieties of~$\P$ in general position
such~that $\codim_{\Bbb{C}}\mu\!=\!d(n\!+\!1)\!-\!1$.
Then the number of degree-$d$ genus-one curves that have 
a fixed generic complex structure on the normalization and 
pass through the constraints~$\mu$ is given~by
$$n_{1,d}(\mu)=\frac{1}{2}\big(
\RT_{1,d}(\mu_1;\mu_2,\ldots,\mu_N)-CR_1(\mu)\big),
\quad\hbox{where}\quad CR_1(\mu)=N(\al_{1,0}),$$
i.e. $CR_1(\mu)$ is the number of zeros of the affine map
$$\psi_{\al_{1,0},\bar{\nu}}\!: E_{1,0}\!=\!L_1\lra\ev^*T\P, 
\quad 
\psi_{\al_{1,0},\bar{\nu}}(\ups)=\bar{\nu}_{\ups}+\al_{1,0}(\ups),$$
over $\bar{\cal V}_1(\mu)$ for a generic section
$\bar{\nu}\!\in\!\Ga\big(\bar{\cal V}_1(\mu);\ev^*T\P\big)$.
\end{prp}

\noindent
Proposition~\ref{g1_prp} is basically 
the main result of the analytic part of~\cite{I}.
The exact statement is not made in~\cite{I}, but 
it can be deduced from the arguments in~\cite{I} by comparing 
with the methods of~\cite{Z2}.\\

\noindent
The general meaning of Proposition~\ref{g1_prp} is the following. 
The number $\RT_{1,d}(\mu_1;\mu_2,\ldots,\mu_N)$
can be viewed as the ``euler class'' of a bundle~$\Ga^{0,1}$
over a closure~$\bar{C}^{\i}$ of the space~$C^{\i}$
of smooth maps from a fixed elliptic curve that pass through
the constraints $\mu_1,\ldots,\mu_N$; see~\cite{LT}.
Then,
\begin{equation}\label{sect_comp_e}
2n_{1,d}(\mu)=\big|\bar{\partial}^{-1}(0)\cap C^{\i}\big|
=\RT_{1,d}(\mu_1;\mu_2,\ldots,\mu_N)-
\sum{\cal C}_{{\cal M}_{\cal T}(\mu)}(\bar{\partial}),
\end{equation}
where $\big\{{\cal M}_{\cal T}(\mu)\big\}$ are complex
finite-dimensional, usually non-compact, manifolds
that stratify
\hbox{$\bar{\partial}^{-1}(0)\cap(\bar{C}^{\i}\!-\!C^{\i})$}.
Equation~\e_ref{sect_comp_e} is an infinite-dimensional analogue
of (2) of Proposition~\e_ref{euler_prp}.
In the finite-dimensional case, 
computation of a contribution to the euler class
from an $s$-regular stratum~${\cal Z}$ of the zero set of section~$s$
reduces to counting
the zeros of a polynomial map between {\it finite}-rank vector bundles
over~$\bar{\cal Z}$, unless $Z$ is $s$-hollow.
The goal in the infinite-dimensional case under consideration 
is a reduction to the same problem
and involves an adoption of the obstruction-bundle
idea of~\cite{T}.
It turns out that ${\cal C}_{{\cal M}_{\cal T}(\mu)}(\bar{\partial})\!=\!0$
for all but one stratum~${\cal M}_{\cal T}(\mu)$
of $\bar{\partial}^{-1}(0)\cap(\bar{C}^{\i}\!-\!C^{\i})$.
The number $CR_1(\mu)$ described by Proposition~\ref{g1_prp}
is the contribution ${\cal C}_{{\cal M}_{\cal T}(\mu)}(\bar{\partial})$
from the only stratum ${\cal M}_{\cal T}(\mu)$ of  
$\bar{\partial}^{-1}(0)\cap(\bar{C}^{\i}\!-\!C^{\i})$
that does contribute
to the ``euler class'' $\RT_{1,d}(\mu_1;\mu_2,\ldots,\mu_N)$
of~$\Ga^{0,1}$.\\

\noindent
As Subsection~\ref{top_sec} suggests, computation of $N(\al_{1,0})$
may require going through a possibly large tree of steps.
We construct this tree as follows.
Each node is a tuple $\si\!=\!(r;k,m;\phi)$, 
where $r\!\ge\!0$ is the distance to the root $\si_0\!=\!(0;1,0;\cdot)$,
$k\!\ge\!1$, and $m\!\ge\!0$.
The tree satisfies the following properties.
If $r\!>\!0$ and $\si^*\!=\!(r\!-\!1;k^*,m^*;\phi^*)$
is the node from which $\si$ is directly descendent, 
we require that $k^*\!\le\!k$, $m^*\!\le\!m$, and at least
one of the inequalities is strict.
Furthermore, $\phi$ specifies a splitting of the set~$[k]$
into $k^*$-disjoint subsets and an assignment of $m\!-\!m^*$
of the elements of the set 
\hbox{$[m]\!=\!\big\{(1,1),\ldots,(1,m)\big\}$} to these subsets.
This description inductively constructs an infinite tree.
However, we will need to consider only the nodes 
$\si\!=\!(r;k,m;\phi)$ \hbox{with $2k\!+\!m\!\le\!n\!+\!1$}.
We will write $\si\!\vdash\!\si'$  to indicate that 
$\si$ is directly descendent from~$\si'$.\\

\noindent
For each node in the above tree, we now define a linear map
between vector bundles over an ms-manifold.
If $\si\!=\!(r;k,m;\phi)$, 
let $\{\si_s\!=\!(s;k_s,m_s;\phi_s)\!:0\!\le\!s\le\!r\}$ 
be the sequence of nodes such that
$\si_r\!=\!\si$ and $\si_s\!\vdash\!\si_{s-1}$ for all $s\!>\!0$.
Put
\begin{gather*}
\bar{\cal V}_{\si}=\bar{\cal V}_{k,m}(\mu),\quad
E_{\si}\!=\!E_{k,m}\lra\bar{\cal V}_{\si},\quad
\al_{\si}=\al_{k,m},\quad
{\cal X}_{\si}={\cal Y}_{\si}\!\times\!\bar{\cal V}_{\si},\quad
{\cal X}_{\si,s}={\cal Y}_{\si,s}\!\times\!\bar{\cal V}_{\si},\\
\quad\hbox{where}\quad
{\cal Y}_{\si}={\cal Y}_{\si,r},\quad
{\cal Y}_{\si,0}=\{pt\},\quad
{\cal Y}_{\si,s}=\Bbb{P}F_{\si_s}\!\times\!{\cal Y}_{\si,s-1}
\hbox{~~if~~}s\!>\!0,\\
\bar{\cal M}_{\si}=
\prod_{i\in\Im~\!\phi}\!\!\bar{\cal M}_{0,i+\phi^{-1}(i)}, \quad
F_{\si}=\bigoplus_{i\in\Im~\!\phi}\!\!\ga_{\si;i}
\lra\bar{\cal M}_{\si}.
\end{gather*}
For the purposes of the last line above,
we view $\phi$ as a map from $[k]\!-\![k^*]$ and 
a subset of~$[m]$ to~$[k^*]$ in the notation of the previous paragraph.
Then, $\ga_{\si;i}\!\lra\!\bar{\cal M}_{0,i+\phi^{-1}(i)}$
is the tautological line bundle; see Subsection~\ref{notation}.
Denote by $\ga_{F_{\si,0}}$ the (trivial) line bundle over~${\cal Y}_{\si,0}$.
Let
$${\cal O}_{\si}={\cal O}_{\si,r},\quad
{\cal O}_{\si,0}=\ev^*T\P,\quad
{\cal O}_{\si,s}={\cal O}_{\si,s-1}\big/\hbox{Im}~\!\bar{\nu}_{\si,s-1}
~~\hbox{if}~s\!>\!0,$$
where
$\bar{\nu}_{\si,s}\in\Ga\big({\cal X}_{\si,s};
\hbox{Hom}(\ga_{F_{\si_s}},{\cal O}_{\si,s})\big)$
is a generic section.
Since $k_{s-1}\!\le\!k_s$, $m_{s-1}\!\le\!m_s$, 
and one of the inequalities is strict,
$$\frac{1}{2}\dim{\cal X}_{\si,s}
\le \frac{1}{2}\dim{\cal X}_{\si}
= \big(n\!+\!1\!-\!2k\!-\!m\big)+
\sum_{s=1}^{s=r}\big(\big|\Im~\!\phi_s\big|\!-\!1\big)
=n-k-r<\rk~\!{\cal O}_{\si,0}-r.$$
Thus, we see inductively that 
each bundle ${\cal O}_{\si,s}$ is well-defined
and a generic section $\bar{\nu}_{\si,s}$ of 
$\hbox{Hom}(\ga_{F_{\si,s}},{\cal O}_{\si,s})$
does not vanish.
Let $\pi_{\si}\!:\ev^*T\P\!\lra\!{\cal O}_{\si}$ be the projection map.
We define 
$$\tilde{\al}_{\si}\in\Ga\big({\cal X}_{\si};
\hbox{Hom}(\ga_{F_{\si}}^*\!\otimes\!E_{\si};
\ga_{F_{\si}}^*\!\otimes\!{\cal O}_{\si})\big),
\quad\hbox{by}\quad
\big\{\tilde{\al}_{\si}(\tau\!\otimes\!\ups)\big\}(w)
=\tau(w)\cdot\pi_{\si}\al_{\si}(\ups)
\in{\cal O}_{\si}.$$
Note that $\tilde{\al}_{\si_0}\!=\!\al_{1,0}$.

\begin{lmm}
\label{comp_top1}
For every node $\si^*$,
$$N\big(\tilde{\al}_{\si^*}\big)=
\big\lan c\big(\ga_{F_{\si^*}}^*\!\otimes\!{\cal O}_{\si^*}\big)
c\big(\ga_{F_{\si^*}}^*\!\otimes\!E_{\si^*}\big)^{-1},
\big[{\cal X}_{\si^*}\big]\big\ran
-\sum_{\si\vdash\si^*}N\big(\tilde{\al}_{\si}\big).$$
\end{lmm}

\noindent
{\it Remark:} For a dense open subset of tuples $\{\bar{\nu}_{\si,s}\}$,
the corresponding linear map~$\al_{\si}$ constructed above is regular
and~$N(\al_{\si})$ is independent of the choice of~$\{\bar{\nu}_{\si,s}\}$.
What we need is that for every bubble type ${\cal T}$ such that
${\cal U}_{\cal T}(\mu)\!\subset\!\bar{\cal V}_{k_r,m_r}(\mu)$
the intersection of the image of the linear map
$$\al_{\cal T}\in\Ga\Big({\cal Y}_{\si}\!\times\!{\cal U}_{\cal T}(\mu);
\hbox{Hom}\big(\!\bigoplus_{i\in\chi({\cal T})}\!\!L_i{\cal T},
\ev^*T\P\big)\Big),
\quad
\al_{\cal T}(\ups)=\!\sum_{i\in\chi({\cal T})}\!\!
{\cal D}_{{\cal T},i}\ups_i,$$
with the subbundle
$$\Im~\!\bar{\nu}_{\si,0}\oplus\ldots\oplus\Im~\!\bar{\nu}_{\si,r-1}
\subset{\cal O}_{\si,0}=\ev^*T\P$$
is $\{0\}$.
The fact that this condition is satisfied for a dense open subset
of tuples~$\{\bar{\nu}_{\si,s}\}$ follows by a dimension count as above,
along with an argument similar to the proof of Lemma~3.10
in~\cite{Z2}.\\

\noindent
{\it Proof of Lemma~\ref{comp_top1}:} (1) By Lemma~\ref{zeros_main},
\begin{equation}\label{comp_top1_e1}
N\big(\tilde{\al}_{\si^*}\big)=
\big\lan c\big(\ga_{F_{\si^*}}^*\!\otimes\!{\cal O}_{\si^*}\big)
c\big(\ga_{F_{\si^*}}^*\!\otimes\!E_{\si^*}\big)^{-1},
\big[{\cal X}_{\si^*}\big]\big\ran
-{\cal C}_{{{\tilde{\al}_{\si^*}}\!'}^{\!-1}(0)}
\big({{\tilde{\al}_{\si^*}}\!'}^{\!\perp}\big).
\end{equation}
Let $\si^*\!=\!(r^*;k^*,m^*;\phi^*)$.
By Lemma~\ref{str_lmm}, 
${{\tilde{\al}_{\si^*}}\!'}^{-1}(0)$ is the union of the sets
$${\cal Z}_{\cal T}^J\equiv{\cal Y}_{\si^*}\times 
\Big(\Bbb{P}E{\cal T}^J-
 \bigcup_{J'\varsubsetneq J}\!\!\Bbb{P}E{\cal T}^{J'}\Big),
\quad\hbox{where}\quad
E^J{\cal T}=\bigoplus_{i\in J}L_i{\cal T}
\lra{\cal U}_{\cal T}(\mu),$$
taken over non-basic bubble types
${\cal T}\!=\!(S^2,[N]\!-\!M_0,I;j,\under{d})$,
with
$|I\!-\!\hat{I}|\!=\!k^*$, $|M_0|\!=\!m^*$, and $\sum d_i\!=\!d$,
and nonempty subsets~$J$ of $I\!-\!\hat{I}\!-\!\chi({\cal T})$.\\
(2) The map $\ga_{\cal T}^{\mu}$ of Proposition~\ref{rat_str_prp}
induces an orientation-preserving homeomorphism~$\ga_{{\cal Z}_{\cal T}}^J$ 
between open neighborhoods of ${\cal Z}_{\cal T}^J$~in 
$${\cal NZ}_{\cal T}^J\equiv {\cal FT}
\oplus \ga_{E{\cal T}^J}^*\!\otimes\!\big(
E{\cal T}^{I-\hat{I}-\chi({\cal T})-J}\!\oplus\! E{\cal T}^{\perp}\big)
\lra {\cal Z}_{\cal T}^J$$
and in ${\cal Y}_{\si^*}\!\times\!\Bbb{P}E_{\si^*}$.
Furthermore, for some $\de,C\!\in\!C^{\i}({\cal Z}_{\cal T}^J;\Bbb{R}^+)$,
with appropriate identifications,
\begin{gather}
\label{comp_top1_e3}
\big|\tilde{\al}_{\si^*}'\big(\ga_{{\cal Z}_{\cal T}}^J(b;\ups,u)\big)-
\al_{{\cal Z}_{\cal T}}^J\big(\rho_{{\cal Z}_{\cal T}}^J(b;\ups,u)\big)\big|
\le C(b)|\ups|^{\frac{1}{p}}
\big|\rho_{{\cal Z}_{\cal T}}^J(b;\ups,u)\big|
\qquad\forall(\ups,u)\!\in\!{\cal NZ}_{{\cal T},\de}^J,\\
\hbox{where}\qquad
\rho_{{\cal Z}_{\cal T}}^J\!:{\cal NZ}_{\cal T}^J\lra
\tilde{\cal N}{\cal Z}_{\cal T}^J\equiv \!
\bigoplus_{h\in\chi({\cal T})}\!\!\tilde{\cal N}_h{\cal Z}_{\cal T}^J,
\notag\\
\tilde{\cal N}_h{\cal Z}_{\cal T}^J=
\begin{cases}
L_{\tilde{\io}_h}^*{\cal T}\!\otimes\!L_h{\cal T},
            &\hbox{if}~h\!\in\!\hat{I}, \tilde{\io}_h\!\in\!J;\\
\ga_{E{\cal T}^J}^*\!\otimes\!L_h{\cal T},&\hbox{otherwise};
\end{cases}\quad
\rho_{{\cal Z}_{\cal T};h}^J(\ups,u)=
\begin{cases}
\rho_{{\cal T};h}(\ups),
&\hbox{if}~h\!\in\!\hat{I}, \tilde{\io}_h\!\in\!J;\\
u_{\tilde{\io}_h}\!\otimes\!\rho_{{\cal T};h}(\ups),
&\hbox{if}~h\!\in\!\hat{I}, \tilde{\io}_h\!\not\in\!J;\\
u_h&\hbox{if}~h\!\in\!\chi({\cal T})\!-\!\hat{I};
\end{cases}\notag
\end{gather}
\begin{gather*}
\al_{{\cal Z}_{\cal T}}^J\in\Ga\Big(
{\cal Z}_{\cal T}^J;\hbox{Hom}\big(\tilde{\cal N}{\cal Z}_{\cal T}^J,
\hbox{Hom}(\ga_{F_{\si^*}}^*\!\otimes\!\ga_{E{\cal T}^J},
\ga_{F_{\si^*}}^*\!\otimes\!{\cal O}_{\si^*})\big)\Big)\notag\\
\begin{split}
\big\{\big\{
\al_{{\cal Z}_{\cal T}}^J\big(\rho_{{\cal Z}_{\cal T}}^J(\ups,u)\big)\big\}
(\tau\!\otimes\!\tilde{\ups})\big\}(w)
=&\tau(w)\cdot\pi_{\si^*}\Big(
\sum_{i\in J}\big\{\al_{{\cal T},i}(\rho_{\cal T}\ups)\big\}(\tilde{\ups}_i)
\\
&+\sum_{i\in\chi({\cal T})-\hat{I}}\!\!\!\!
{\cal D}_{{\cal T},i}(u_i\tilde{\ups}) 
+\sum_{i\in I-\hat{I}-\chi({\cal T})-J}\!\!\!\!\!\!\!
\big\{\al_{{\cal T},i}(\rho_{\cal T}\ups)\big\}(u_i\tilde{\ups})\Big)
\in {\cal O}_{\si^*}.
\end{split}\end{gather*}
Above $\rho_{{\cal T};h}$ denotes the $h$th component of $\rho_{\cal T}$,
i.e.~$\tilde{\ups}_h$ in the notation of Proposition~\ref{rat_str_prp}.
Note that the section $\al_{{\cal Z}_{\cal T}}^J$ is well-defined.
By Lemma~\ref{str_lmm} and the splitting~\e_ref{cart_split1a},
possibly applied several times,
$\pi_{\bar{\nu}_{\si^*}}^{\perp}\!\circ\al_{{\cal Z}_{\cal T}}^J$
has full rank on every fiber of $\tilde{\cal N}{\cal Z}_{\cal T}^J$,
provided the sections $\big\{\bar{\nu}_{\si,s}\!:0\!\le\! s\!\le\!r^*\big\}$
are generic.
Then by~\e_ref{comp_top1_e3},
\begin{equation}\label{comp_top1_e5}
\big|\pi_{\bar{\nu}_{\si^*}}^{\perp}\!\circ\tilde{\al}_{\si^*}'
\big(\ga_{{\cal Z}_{\cal T}}^J(b;\ups,u)\big)-
\pi_{\bar{\nu}_{\si^*}}^{\perp}\!\circ\al_{{\cal Z}_{\cal T}}^J
\big(\rho_{{\cal Z}_{\cal T}}^J(b;\ups,u)\big)\big|
\le C(b)|\ups|^{\frac{1}{p}}\big|\pi_{\bar{\nu}_{\si^*}}^{\perp}
\!\circ\al_{{\cal Z}_{\cal T}}^J
\big(\rho_{{\cal Z}_{\cal T}}^J(b;\ups,u)\big)\big|.
\end{equation}
for all $(\ups,u)\!\in\!{\cal NZ}_{{\cal T},\de}^J$.
Thus, $\pi_{\bar{\nu}_{\si^*}}^{\perp}\!\circ\al_{{\cal Z}_{\cal T}}^J$ is
the resolvent for 
$\ga_{{\cal Z}_{\cal T}}^{J*}
 (\pi_{\bar{\nu}_{\si^*}}^{\perp}\!\circ\tilde{\al}_{\si^*}')$;
see Definition~\ref{top_dfn1b}.
If ${\cal T}$ is not semiprimitive or 
$J\!\neq\!I\!-\!\hat{I}\!-\!\chi({\cal T})$,
the rank of $\tilde{\cal N}{\cal Z}_{\cal T}^J$
is less than the rank of~${\cal NZ}_{\cal T}^J$.
It follows that ${\cal Z}_{\cal T}^J$ is 
$\pi_{\bar{\nu}_{\si^*}}^{\perp}\!\circ\tilde{\al}_{\si^*}'$-hollow and
\begin{equation}\label{comp_top1_e7a}
{\cal C}_{{\cal Z}_{\cal T}^J}
\big({{\tilde{\al}_{\si^*}}\!\!'}^{\!\perp}\big)=0
\qquad\hbox{if}~~
{\cal T} \hbox{ is not semiprimitive or }
J\neq I\!-\!\hat{I}\!-\!\chi({\cal T})
\end{equation}
by Proposition~\ref{euler_prp}.\\
(3)  On the other hand,
by the analytic estimate~\e_ref{comp_top1_e5}, Proposition~\ref{euler_prp},
and the splitting~\e_ref{cart_split1b}, 
\begin{gather}\label{comp_top1_e7b}
{\cal C}_{{\cal Z}_{\cal T}^J}
\big({{\tilde{\al}_{\si^*}}\!\!'}^{\!\perp}\big)=
N\big(\al_{\si^*,{\cal T}}\big)
\qquad\hbox{if}~~
{\cal T} \hbox{ is semiprimitive and }
J= I\!-\!\hat{I}\!-\!\chi({\cal T}),\\
\hbox{where}\quad
\al_{\si^*,{\cal T}}\in\Ga\big(
{\cal Y}_{\si^*,{\cal T}}\!\times\!\bar{\cal U}_{\bar{\cal T}}(\mu);
\hbox{Hom}({\cal NZ}_{{\si}^*,\cal T},
\ga_{E{\cal T}}^*\!\otimes\!{\cal O}_{\si^*,{\cal T}})\big),\quad
{\cal Y}_{\si^*,{\cal T}}\!=\!{\cal Y}_{\si^*}\!\times\!\Bbb{P}E{\cal T},
\notag
\end{gather} 
\begin{gather}
E{\cal T}\!=\!\bigoplus_{i\in\hat{I}}\ga_{{\cal T};i}
\lra {\cal M}_{\si^*,{\cal T}}\!=\!\!
\prod_{i\in I-\chi({\cal T})}\!\!\!\!\!\!
\bar{\cal M}_{0,H_i{\cal T}+M_i{\cal T}},
\quad
{\cal NZ}_{{\si}^*,\cal T}=
\bigoplus_{h\in\hat{I}}\ga_{{\cal T};\io_h}^*\!\otimes\!L_h\bar{\cal T} 
\oplus\!\bigoplus_{h\in\chi({\cal T})-\hat{I}}\!\!\!\!\!\!
\ga_{E{\cal T}}^*\!\otimes\!L_h\bar{\cal T},\notag\\
{\cal O}_{\si^*,{\cal T}}=
{\cal O}_{\si^*}\big/\hbox{Im}~\!\bar{\nu}_{\si^*}
\approx \ga_{F_{\si^*}}\!\otimes\!\big((\ga_{F_{\si^*}}^*
\!\otimes\!{\cal O}_{\si^*})\big/\Bbb{C}\bar{\nu}_{\si^*}\big),\notag\\
\big\{\al_{\si^*,{\cal T}}(u\!\otimes\!\ups)\big\}(\tilde{\ups})=
\pi_{\si^*,{\cal T}}
\Big(\sum_{i\in\hat{I}}u_i(\tilde{\ups}_i)
 \big({\cal D}_{\bar{\cal T},i}\ups_i\big)
+\sum_{i\in\chi({\cal T})-\hat{I}}\!\!\!\!
u_i(\tilde{\ups})\big({\cal D}_{\bar{\cal T},i}\ups_i\big)\Big)
\in{\cal O}_{\si^*,{\cal T}}, \notag
\end{gather}
$\ga_{{\cal T};i}\lra\bar{\cal M}_{0,H_i{\cal T}+M_i{\cal T}}$
is the tautological line bundle, and
\hbox{$\pi_{\si^*,{\cal T}}\!:\ev^*T\P\!\lra\!{\cal O}_{\si^*,{\cal T}}$}
is the quotient projection map.
We next observe that
\begin{gather}\label{comp_top1_e9}
N\big(\al_{\si^*,{\cal T}}\big)=
N\big(\tilde{\al}_{\si^*,{\cal T}}\big),\qquad\hbox{where}\\
\tilde{\al}_{\si^*,{\cal T}}\in\Ga\big(
{\cal Y}_{\si^*,{\cal T}}\!\times\!\bar{\cal U}_{\bar{\cal T}}(\mu);
\hbox{Hom}(\ga_{E{\cal T}}^*\!\otimes\! E_{k,m},
\ga_{E{\cal T}}^*\!\otimes\!{\cal O}_{\si^*,{\cal T}})\big),\notag\\
k=|\chi({\cal T})|=|\bar{I}|,\quad 
m=m^*+\!\sum_{i\in I-\chi({\cal T})}\!\!\!\!|M_i{\cal T}|
\quad\hbox{and}\quad
\big\{\tilde{\al}_{\si^*,{\cal T}}(\tau\!\otimes\!\ups)\big\}(w)
=\tau(w)\cdot\pi_{\si^*,{\cal T}}\al_{k,m}(\ups).\notag
\end{gather}
The reason for the equality~\e_ref{comp_top1_e9} is the following.
For a generic $\bar{\nu}\!\in\!\Ga\big(
{\cal Y}_{\si^*,{\cal T}}\!\times\!\bar{\cal U}_{\bar{\cal T}}(\mu);
\ga_{E{\cal T}}^*\!\otimes\!{\cal O}_{\si^*,{\cal T}}\big)$,
the affine~maps
$$\psi_{\al_{\si^*,{\cal T}},\bar{\nu}}
\!\equiv\!\bar{\nu}\!+\!\al_{\si^*,{\cal T}}\!:
{\cal NZ}_{{\si}^*,\cal T}\!\lra\!
\ga_{E{\cal T}}^*\!\otimes\!{\cal O}_{\si^*,{\cal T}}
~~~\hbox{and}~~~
\psi_{\tilde{\al}_{\si^*,{\cal T}},\bar{\nu}}
\!\equiv\!\bar{\nu}\!+\!\tilde{\al}_{\si^*,{\cal T}}\!:
\ga_{E{\cal T}}^*\!\otimes\! E_{k,m}\!\lra\!
\ga_{E{\cal T}}^*\!\otimes\!{\cal O}_{\si^*,{\cal T}}$$
have no zeros over the complement of ${\cal Z}_{\cal T}^J$,
since it is a finite union of smooth manifolds of dimension
less than that of~${\cal Z}_{\cal T}^J$.
There is a canonical identification of the line bundle $\ga_{E{\cal T}}^*$
with each line bundle~$\ga_{{\cal T};i}$ over~${\cal Z}_{\cal T}^J$.
This identification induces a bijection between 
the zeros of the two affine maps that lie over~${\cal Z}_{\cal T}^J$.
The identity~\e_ref{comp_top1_e9} follows from this argument along~with
$$N\big(\al_{\si^*,{\cal T}}\big)=
^{\pm}\!\big|\psi_{\al_{\si^*,{\cal T}},\bar{\nu}}^{-1}(0)\big|
\quad\hbox{and}\quad
N\big(\tilde{\al}_{\si^*,{\cal T}}\big)=
^{\pm}\!\big|\psi_{\tilde{\al}_{\si^*,{\cal T}},\bar{\nu}}^{-1}(0)\big|.$$
(4) From equations~\e_ref{comp_top1_e1}, 
\e_ref{comp_top1_e7a}, \e_ref{comp_top1_e7b}, and \e_ref{comp_top1_e9},
we conclude that that
\begin{equation*}\begin{split}
N\big(\tilde{\al}_{\si^*}\big)&=
\big\lan c\big(\ga_{F_{\si^*}}^*\!\otimes\!{\cal O}_{\si^*}\big)
c\big(\ga_{F_{\si^*}}^*\!\otimes\!E_{\si^*}\big)^{-1},
\big[{\cal X}_{\si^*}\big]\big\ran
-\!\!\sum_{(k,m)>(k^*,m^*)}
\sum_{~|\chi({\cal T})|=k,\!
\sum\limits_{i\in I-\chi({\cal T})}\!\!\!\!\!|M_i{\cal T}|=m-m^*}
\!\!\!\!\!\!\!\!\!\!\!\!N\big(\tilde{\al}_{\si^*,{\cal T}}\big)\notag\\
&=\big\lan c\big(\ga_{F_{\si^*}}^*\!\otimes\!{\cal O}_{\si^*}\big)
c\big(\ga_{F_{\si^*}}^*\!\otimes\!E_{\si^*}\big)^{-1},
\big[{\cal X}_{\si^*}\big]\big\ran
-\!\!\sum_{\si\vdash\si^*}N\big(\tilde{\al}_{\si}\big).
\end{split}\end{equation*}
The inner sum on the first line above is taken over all equivalence 
classes of semiprimitive bubble types 
\hbox{${\cal T}\!=\!(S^2,N\!-\!M_0,I;j,\under{d})$}
such that $|I\!-\!\hat{I}|\!=\!k^*$,  $|M_0|\!=\!m^*$,
and $\sum d_i\!=\!d$.
Condition \hbox{$(k,m)\!>\!(k^*,m^*)$} means that $k\!\ge\!k$,
$m\!\ge\!m^*$ and one of the inequalities is strict.

\begin{lmm}
\label{comp_top2}
For every node $\si\!=\!(r;k,m;\phi)$ and positive integer $s\!\le\!r\!-\!1$,
$$\big\lan c\big({\cal O}_{\si,s+1}\big)c\big(E_{\si}\big)^{-1},
\big[{\cal X}_{\si,s}\big]\big\ran
=\big\lan c\big({\cal O}_{\si,s}\big)c\big(E_{\si}\big)^{-1},
\big[{\cal X}_{\si,s-1}\big]\big\ran,$$
where $\{\si_s\}$ is the sequence corresponding to $\si$
defined in the paragraph preceding Lemma~\ref{comp_top1}.
\end{lmm}

\noindent
{\it Proof:} Since 
${\cal O}_{\si,s+1}\!\approx\!{\cal O}_{\si,s}\big/\ga_{F_{\si_s}}$,
\begin{equation}\label{comp_top2e1}
\big\{c\big({\cal O}_{\si,s+1}\big)c\big(E_{\si}\big)^{-1}
\big\}_{\dim {\cal X}_{\si,s}}=
\sum_{l=0}^{\dim{\cal X}_{\si,s}}
\!\sum_{l_1+l_2=l}\!\la_{F_{\si_s}}^{l_1}c_{l_2}\big({\cal O}_{\si,s}\big)
\big\{c\big(E_{\si}\big)^{-1}\big\}_{\dim{\cal X}_{\si,s}-2l}.
\end{equation}
By construction, $\la_{F_{\si_s}}\!\in\!H^*(\Bbb{P}F_{\si_s})$,
while 
$c({\cal O}_{\si,s}),c(E_{\si})\!\in\!H^*({\cal X}_{\si,s-1})$.
Thus, \e_ref{comp_top2e1} gives
\begin{equation}\label{comp_top2e3}\begin{split}
\big\{c\big({\cal O}_{\si,s+1}\big)c\big(E_{\si}\big)^{-1}
\big\}_{\dim {\cal X}_{\si,s}}
&=\la_{F_{\si_s}}^{n_{\si}}\sum_{l=0}^{\dim{\cal X}_{\si,s}}
\! c_{l-n_{\si}}\big({\cal O}_{\si,s}\big)
\big\{c\big(E_{\si}\big)^{-1}\big\}_{\dim{\cal X}_{\si,s}-2l}\\
&=\la_{F_{\si_s}}^{n_{\si}}
\big\{c\big({\cal O}_{\si,s}\big)c\big(E_{\si}\big)^{-1}
\big\}_{\dim{\cal X}_{\si,s-1}},
\end{split}\end{equation}
where $n_{\si}\!=\!\dim\Bbb{P}F_{\si_s}$.
By~\e_ref{zeros_main_e},
\begin{equation}\label{comp_top2e5}
\big\lan\la_{F_{\si_s}}^{n_{\si}},\big[\Bbb{P}F_{\si_s}\big]\big\ran
=\big\lan c\big(F_{\si_s}\big)^{-1},\big[\bar{\cal M}_{\si_s}\big]
\big\ran
=\prod_{i\in\Im~\!\phi_s}
\big\lan c\big(\ga_{\si_s;i}\big)^{-1},
\big[\bar{\cal M}_{0,i+\phi_s^{-1}(i)}\big]\big\ran=1.
\end{equation}
The last identity is a consequence of (1) of Lemma~\ref{ag_lmm}.
The claim follows from~\e_ref{comp_top2e1}-\e_ref{comp_top2e5}.

\begin{crl}
\label{comp_top3}
For every node $\si\!=\!(r;k,m;\phi)$,
$$\big\lan c\big(\ga_{F_{\si}}^*\!\otimes\!{\cal O}_{\si}\big)
c\big(\ga_{F_{\si}}^*\!\otimes\!E_{\si}\big)^{-1},
\big[{\cal X}_{\si}\big]\big\ran
=\big\lan c(\ev^*T\P)c(E_{k,m})^{-1},
\big[\bar{\cal V}_{k,m}(\mu)\big]\big\ran.$$
\end{crl}

\noindent
{\it Proof:}
Since 
$\rk~\!{\cal O}_{\si}\!=\!\rk~\!E_{\si}\!+\!\frac{1}{2}\dim{\cal X}_{\si}$,
we can identify $E_{\si}$ with a subbundle of~${\cal O}_{\si}$.
Then,
\begin{gather}
c\big(\ga_{F_{\si}}^*\!\otimes\!{\cal O}_{\si}\big)
c\big(\ga_{F_{\si}}^*\!\otimes\!E_{\si}\big)^{-1}
=c\big(\ga_{F_{\si}}^*\!\otimes\!{\cal O}_{\si}\big/
\ga_{F_{\si}}^*\!\otimes\!E_{\si}\big)
=c\big(\ga_{F_{\si}}^*\!\otimes\!({\cal O}_{\si}/E_{\si})\big)\Lra\notag\\
\label{comp_top3e1}
\big\{c\big(\ga_{F_{\si}}^*\!\otimes\!{\cal O}_{\si}\big)
c\big(\ga_{F_{\si}}^*\!\otimes\!E_{\si}\big)^{-1}
\big\}_{\dim{\cal X}_{\si}}
=\sum_{l=0}^{\dim{\cal X}_{\si}}\! \la_{F_{\si}}^l
\big\{c({\cal O}_{\si})c(E_{\si})^{-1}\big\}_{\dim{\cal X}_{\si}-2l}.
\end{gather}
Similarly to the proof of Lemma~\ref{comp_top2},
\e_ref{comp_top3e1} gives
\begin{equation}\label{comp_top3e3}\begin{split}
\big\lan c\big(\ga_{F_{\si}}^*\!\otimes\!{\cal O}_{\si}\big)
c\big(\ga_{F_{\si}}^*\!\otimes\!E_{\si}\big)^{-1},
\big[{\cal X}_{\si}\big]\big\ran
&=\big\lan c\big({\cal O}_{\si}\big)
c\big(E_{\si}\big)^{-1},\big[{\cal X}_{\si,r-1}\big]\big\ran\\
&=\big\lan c\big({\cal O}_{\si,r}\big)
c\big(E_{\si}\big)^{-1},\big[{\cal X}_{\si,r-1}\big]\big\ran.
\end{split}\end{equation}
Applying Lemma~\ref{comp_top2} to the last expression in~\e_ref{comp_top3e3} 
and using \hbox{${\cal O}_{\si,1}\!\approx\!(\ev^*T\P)/\Bbb{C}$}, we~obtain
$$\big\lan c\big(\ga_{F_{\si}}^*\!\otimes\!{\cal O}_{\si}\big)
c\big(\ga_{F_{\si}}^*\!\otimes\!E_{\si}\big)^{-1}\!,
\!\big[{\cal X}_{\si}\big]\big\ran
=\big\lan c\big({\cal O}_{\si,1}\big)
c\big(E_{\si}\big)^{-1}\!,\!\big[{\cal X}_{\si,0}\big]\big\ran
=\big\lan c\big(\ev^*T\P)c(E_{k,m})^{-1}\!,
\!\big[\bar{\cal V}_{k,m}(\mu)\big]\big\ran.$$
\\

\noindent
We now combine Lemma~\ref{comp_top1} and Corollary~\ref{comp_top3}
to obtain a topological formula for the number~$N(\al_{1,0})$. 
For any integers~$k$ and~$k^*$, let $\th_{k^*}^k$ denote the number
of ways of splitting a set of $k^*$-elements into $k$ nonempty subsets.
For every pair $(k^*,m^*)\!\ge\!(1,0)$ of integers, 
we define $\Th(k^*,m^*)$ inductively~by
\begin{equation}\label{comp_top_dfn}
\Th(1,0)\!=\!1,\quad
\Th(k^*,m^*)=-\!\!\!\sum_{(1,0)\le(k,m)<(k^*,m^*)}\!\!\!
\binom{m^*}{m}k^{m^*-m}\th_{k^*}^k\Th(k,m)
\hbox{~~if~~}(k^*,m^*)\!>\!(1,0).
\end{equation}

\begin{crl}
\label{comp_top4}
With notation as above,
$$N(\al_{1,0})= \sum_{(1,0)\le(k,m)}\!\!\Th(k,m)
\!\sum_{l=0}^{n+1-(2k+m)}\!\!
\binom{n\!+\!1}{l} \big\lan a^l\tilde{\eta}_{n+1-(2k+m)-l},
\big[\bar{\cal V}_{k,m}(\mu)\big]\big\ran.$$
\end{crl}

\noindent
{\it Proof:} 
By Lemma~\ref{comp_top1} and Corollary~\ref{comp_top3},
\begin{equation}\label{comp_top4e1}
N\big(\al_{1,0}\big)=N\big(\tilde{\al}_{\si_0}\big)
=\sum_{(1,0)\le(k,m)}\!\!\!\!\Th(k,m)
\big\lan c(\ev^*T\P)c(E_{k,m})^{-1},
\big[\bar{\cal V}_{k,m}(\mu)\big]\big\ran.
\end{equation}
Since $E_{k,m}=\bigoplus L_i$,
\begin{equation}\label{comp_top4e2}
c(E_{k,m})^{-1}
=\prod_{i=1}^{i=k}\big(1+c_1(L_i)\big)^{-1}
=\prod_{i=1}^{i=k}\sum_{l=0}^{\i}c_1^l(L_i^*)
=\sum_{l=0}^{\i}\tilde{\eta}_l.
\end{equation}
The last equality above is immediate from the definition of~$\tilde{\eta}_l$;
see Subsection~\ref{notation}.
The claim follows from~\e_ref{comp_top4e1} and~\e_ref{comp_top4e2},
along with $c(\ev^*T\P)\!=\!(1\!+\!a)^{n+1}$.

\subsection{Combinatorics}

\noindent
In this subsection, we show that the topological expression 
for~$N(\al_{1,0})$ given in Corollary~\ref{comp_top4} is the same
as the topological expression for~$CR_1(\mu)$
given in Proposition~\ref{cr_str}.
This fact is immediate from Corollary~\ref{comp_com4}.
We start by proving an explicit formula for the numbers~$\Th(k,m)$.

\begin{lmm}
\label{comp_com1}
If $(k,m)\!\ge\!(1,0)$, $\Th(k,m)=(-1)^{k+m-1}k^m(k\!-\!1)!$.
\end{lmm}

\noindent
(1) We first start verify this formula in the case $k=1$.
By~\e_ref{comp_top_dfn},
\begin{equation}\label{comp_com1e1}
\Th(1,0)=1,\qquad
\Th(1,m^*)=-\sum_{m=0}^{m^*-1}\binom{m^*}{m}\Th(1,m)
\quad\hbox{if~~}m^*\!>\!(1,0).
\end{equation}
We need to show that $\Th(1,m)\!=\!(-1)^m$.
If $m\!=\!0$, this is the case.
Suppose $m^*\!\ge\!1$ and \hbox{$\Th(1,m)\!=\!(-1)^m$} for all $m\!<\!m^*$.
Then, by~\e_ref{comp_com1e1},
$$\Th(1,m^*)=-\!\sum_{m=0}^{m^*-1}\!\!\binom{m^*}{m}\Th(1,m)
=-\sum_{m=0}^{m^*}\!\binom{m^*}{m}(-1)^m+(-1)^{m^*}
=-(1\!-\!1)^{m^*}\!+\!(-1)^{m^*}\!=\!(-1)^{m^*}\!\!,$$
as needed.\\
(2) We now verify the formula in the general case.
It is easy to see from the definition of $\th_{k^*}^k$
in the previous subsection that
\begin{equation}\label{comp_com1e3}
\th_k^k=1\hbox{~~if~}k\!\ge\!1
\quad\hbox{and}\quad
\th_{k^*}^k=k\th_{k^*-1}^k+\th_{k^*-1}^{k-1}
\hbox{~~if~}k\!\ge\!2.
\end{equation}
Suppose $k^*\!\ge\!2$ and the claimed formula holds
for all $(k,m)$ with $(1,0)\!\le\!(k,m)\!<\!(k^*,m^*)$.
Then by~\e_ref{comp_top_dfn},
\begin{equation}\label{comp_com1e5}\begin{split}
\Th(k^*,m^*)
&=-\!\!\sum_{(1,0)\le(k,m)<(k^*,m^*)}\!\!\!
\binom{m^*}{m}k^{m^*-m}\th_{k^*}^k\Th(k,m)\\
&=k^{m^*}\!\!\!\!\sum_{(1,0)\le(k,m)<(k^*,m^*)}\!\!\!\!\!\!
(-1)^{k+m}\binom{m^*}{m}\th_{k^*}^k(k\!-\!1)!
\end{split}\end{equation}
Using~\e_ref{comp_com1e3}, we obtain
\begin{alignat}{1}
&\sum_{(1,0)\le(k,m)<(k^*,m^*)}\!\!\!\!\!\!\!\!\!
(-1)^{k+m}\binom{m^*}{m}\th_{k^*}^k(k\!-\!1)!
=\!\!\sum_{(1,0)\le(k,m)<(k^*,m^*)}\!\!\!\!\!\!\!\!\!
(-1)^{k+m}\binom{m^*}{m}
\big(k\th_{k^*-1}^k+\th_{k^*-1}^{k-1}\big)(k\!-\!1)!\notag\\
\label{comp_com1e7}
&\qquad\qquad\qquad\qquad\qquad\qquad\qquad
=\!\sum_{m=0}^{m^*-1}\!(-1)^m\binom{m^*}{m}
\sum_{k=1}^{k^*}(-1)^k\big(\th_{k^*-1}^kk!+\th_{k^*-1}^{k-1}(k\!-\!1)!\big)\\
&\qquad\qquad\qquad\qquad\qquad\qquad\qquad\qquad
+(-1)^{m^*}\!\sum_{k=1}^{k^*-1}\!(-1)^k
\big(\th_{k^*-1}^kk!+\th_{k^*-1}^{k-1}(k\!-\!1)!\big).\notag
\end{alignat}
Note that
\begin{alignat}{1}
\label{comp_com1e9}
\sum_{k=1}^{k^*}(-1)^k\big(\th_{k^*-1}^kk!+\th_{k^*-1}^{k-1}(k\!-\!1)!\big)
&=\sum_{k=1}^{k^*}(-1)^k\th_{k^*-1}^kk!
-\sum_{k=0}^{k^*-1}(-1)^k\th_{k^*-1}^kk!
=0;\\
\sum_{k=1}^{k^*-1}(-1)^k
\big(\th_{k^*-1}^kk!+\th_{k^*-1}^{k-1}(k\!-\!1)!\big)
&=\sum_{k=1}^{k^*-1}(-1)^k\th_{k^*-1}^kk!
-\sum_{k=0}^{k^*-2}(-1)^k\th_{k^*-1}^kk!
=(-1)^{k^*-1}(k^*\!-\!1)!,\notag
\end{alignat}
since $c_{k^*-1}^{k^*}\!=\!0$, $c_{k^*-1}^{k^*-1}\!=\!1$,
and $c_{k^*-1}^0\!=\!0$ if $k^*\!>\!1$.
Combining equations~\e_ref{comp_com1e5}-\e_ref{comp_com1e9},
we verify the claimed identity for $(k,m)\!=\!(k^*,m^*)$.\\

\noindent
We next need to relate the intersection numbers 
$a^l\tilde{\eta}_{l'}$ and~$a^l\tilde{\eta}_{l'}$.
We break the computation into several steps.

\begin{lmm}
\label{comp_com2}
Suppose ${\cal T}\!=\!\big(S^2,M,I;j,\under{d})$
is a basic bubble type, $i\!\in\!I$, and $M_i\!\subset\!M_i{\cal T}$.
Then, under the splitting~\e_ref{cart_split2},
with $\bar{\cal T}\!=\!{\cal T}/M_i$,
$$c_1(L_{i'}^*{\cal T})\big|\bar{\cal U}_{{\cal T}(M_i)}(\mu)
=\begin{cases}
\ga_{{\cal T};i}^*\!\times\!1,&\hbox{if~}i'\!=\!i;\\
1\!\times\!c_1(L_{i'}^*\bar{\cal T}),&\hbox{if~}i'\!\neq\!i;
\end{cases}\qquad
c_1({\cal L}_{i'}^*{\cal T})\big|\bar{\cal U}_{{\cal T}(M_i)}(\mu)
=1\!\times\!c_1({\cal L}_{i'}^*\bar{\cal T}).$$
\end{lmm}

\noindent
{\it Proof:} The first identity and the case $i'\!\neq\!i$
of the second identity are immediate from the definitions.
In the remaining case, by~\e_ref{normal_bundle1}, we~have
\begin{equation}\label{comp_com2e1}
c_1({\cal L}_i^*{\cal T})\big|\bar{\cal U}_{{\cal T}(M_i)}(\mu)
=c_1(L_i^*{\cal T})\big|\bar{\cal U}_{{\cal T}(M_i)}(\mu)
-\sum_{\eset\neq M_i'\subset M_i{\cal T}}\!\!\!\!\!
\hbox{PD}_{\bar{\cal U}_{\cal T}(\mu)}
\bar{\cal U}_{{\cal T}(M_i')}(\mu)\big|
\bar{\cal U}_{{\cal T}(M_i)}(\mu).
\end{equation}
By definition of the spaces, 
\begin{equation}\label{comp_com2e3}
\hbox{PD}_{\bar{\cal U}_{\cal T}(\mu)}
\bar{\cal U}_{{\cal T}(M_i')}(\mu)\big|
\bar{\cal U}_{{\cal T}(M_i)}(\mu)=
\begin{cases}
0,&\hbox{if}~M_i'\!\not\subset\! M_i\hbox{~and~}M_i\!\not\subset\!M_i';\\
1\!\times\!\hbox{PD}_{\bar{\cal U}_{\bar{\cal T}}(\mu)}      
\bar{\cal U}_{\bar{\cal T}(M_i'-M_i)}(\mu),
&\hbox{if}~M_i\!\varsubsetneq\! M_i';\\
\hbox{PD}_{\bar{\cal U}_{{\cal T}_0}}
\bar{\cal U}_{{\cal T}_0(M_i-M_i')}\!\times\!1
,&\hbox{if}~M_i'\!\varsubsetneq\! M_i.
\end{cases}
\end{equation}
where ${\cal T}_0\!=\!(S^2,\hat{1}\!+\!M_i,\{i\};i,0)$,
i.e.~$\bar{\cal U}_{{\cal T}_0}\!=\!\bar{\cal M}_{0,\hat{1}+M_i}$.
Plugging~\e_ref{comp_com2e3}, (2) of Lemma~\ref{ag_lmm},
and the case $i'\!=\!i$ of the first statement of
this lemma into~\e_ref{comp_com2e1},
we obtain the remaining claim.

\begin{crl}
\label{comp_com3}
For all $k\!\ge\!1$, $m\!\ge\!0$, and $l\!\ge\!0$,
$$\big\lan a^l\tilde{\eta}_{n+1-(2k+m)-l},
                      \big[\bar{\cal V}_{k,m}(\mu)\big]\big\ran
=\!\sum_{m^*\ge m}\!\!\binom{m^*}{m}k^{m^*-m}
\big\lan a^l\eta_{n+1-(2k+m^*)-l},
\big[\bar{\cal V}_{k,m^*}(\mu)\big]\big\ran.$$
\end{crl}

\noindent
{\it Proof:} Let ${\cal T}\!=\!(S^2,[N]\!-\!M_0,I;j,\under{d})$
be a basic bubble type such that $|I|\!=\!k$, $|M_0|\!=\!m$,
and $\sum d_i\!=\!d$.
By Lemma~\ref{comp_com2} and (1) of Lemma~\ref{ag_lmm},
\begin{equation}\label{comp_com3e1}
\big\lan a^l\tilde{\eta}_{n+1-(2k+m)-l},
                      \big[\bar{\cal U}_{\cal T}(\mu)\big]\big\ran
=\sum_{M_0\subset M_0^*\subset[N]}\!\!\!\!\!
\big\lan a^l\eta_{n+1-(2k+|M_0^*|)-l},
                      \big[\bar{\cal U}_{{\cal T}/M_0^*}(\mu)\big]\big\ran,
\end{equation} 
where ${\cal T}/M_0^*\!=\!(S^2,[N]\!-\!M_0^*,I;j,\under{d})$.	
The claim is obtained by summing~\e_ref{comp_com3e1}
over all equivalence classes of bubble types~${\cal T}$
of the above~form.

\begin{crl}
\label{comp_com4}
For all $k\!\ge\!1$ and $l\!\ge\!0$,
$$\sum_{m\ge0}\Th(k,m)
\big\lan a^l\tilde{\eta}_{n+1-(2k+m)-l},
                      \big[\bar{\cal V}_{k,m}(\mu)\big]\big\ran
=(-1)^{k-1}(k\!-\!1)!\big\lan a^l\eta_{n+1-2k-l},
\big[\bar{\cal V}_k(\mu)\big]\big\ran.$$
\end{crl}

\noindent
{\it Proof:}
By Lemma~\ref{comp_com1} and Corollary~\ref{comp_com3},
\begin{equation*}\begin{split}
&\sum_{m\ge0}\Th(k,m)
\big\lan a^l\tilde{\eta}_{n+1-(2k+m)-l},
                      \big[\bar{\cal V}_{k,m}(\mu)\big]\big\ran\\
&\qquad\qquad\qquad=(-1)^{k-1}(k\!-\!1)!\sum_{m\ge0}\sum_{m^*\ge m}
(-1)^m\binom{m^*}{m}k^{m^*} \big\lan a^l\eta_{n+1-(2k+m^*)-l},
\big[\bar{\cal V}_{k,m^*}(\mu)\big]\big\ran\\
&\qquad\qquad\qquad=(-1)^{k-1}(k\!-\!1)!\sum_{m^*\ge0} k^{m^*}
\bigg(\sum_{m\le m^*}(-1)^m\binom{m^*}{m}\bigg)
 \big\lan a^l\eta_{n+1-(2k+m^*)-l},
\big[\bar{\cal V}_{k,m^*}(\mu)\big]\big\ran\\
&\qquad\qquad\qquad=(-1)^{k-1}(k\!-\!1)!
\big\lan a^l\eta_{n+1-2k-l},
\big[\bar{\cal V}_{k,0}(\mu)\big]\big\ran,
\end{split}\end{equation*}
since
$$\sum_{m\le m^*}(-1)^m\binom{m^*}{m}k^{m^*}
=(1-1)^{m^*}=0\qquad\hbox{if~~}m^*\neq0.$$

\begin{lmm}
\label{ag_lmm}
(1) If $J$ is a finite set of cardinality at least two,
$\big\lan c_1^{|J|-2}(\ga_J^*),\big[\bar{\cal M}_{0,J}\big]
\big\ran\!=\!1$, where \hbox{$\ga_J\!\lra\!\bar{\cal M}_{0,J}$} 
is the tautological line bundle.\\
(2) If ${\cal T}\!=\!(S^2,M,I;j,\under{d})$ is a basic bubble type,
\hbox{$i\!\in\!I$}, and $M_i$ is nonempty subset of~$M_i{\cal T}$,
under the splitting~\e_ref{cart_split2},
$$\hbox{PD}_{\bar{\cal U}_{\cal T}(\mu)}
\bar{\cal U}_{{\cal T}(M_i)}(\mu)\big|
\bar{\cal U}_{{\cal T}(M_i)}(\mu)=
-1\!\times\!c_1(L_i^*\bar{\cal T})
+c_1(\ga_{{\cal T};i}^*)\!\times\!1
-\sum_{\eset\neq M_i'\varsubsetneq M_i}\!\!\!\!\!
\hbox{PD}_{\bar{\cal U}_{{\cal T}_0}}
\bar{\cal U}_{{\cal T}_0(M_i-M_i')}\!\times\!1,$$
where ${\cal T}_0\!=\!(S^2,\hat{1}\!+\!M_i,\{i\};i,0)$
and $\bar{\cal T}\!=\!{\cal T}/M_i$.
\end{lmm}

\noindent
{\it Proof:} (1) Both statements are straightforward consequences
of well-known facts in algebraic geometry; see~\cite{P2}.
In our notation, $\bar{\cal M}_{0,J}$ is the Deligne-Mumford
moduli space of rational curves with points marked by the set
$\{\hat{0}\}\!+\!J$ and $c_1(\ga_J^*)\!=\!\psi_{\hat{0}}$.
Thus, if $j_1,j_2\!\in\!J$ \hbox{and $j_1\!\neq\!j_2$},
\begin{equation}\label{ag_lmm_e1}
c_1(\ga_J^*)=\psi_{\hat{0}}=\!
\sum_{\eset\neq J'\subset J-\{j_1,j_2\}}\!\!\!\!\!\!\!
\hbox{PD}_{\bar{\cal U}_{{\cal T}_0}}
\bar{\cal U}_{{\cal T}_0(J')},
\end{equation}
where ${\cal T}_0\!=\!(S^2,J,\{i\};i,0)$.
Since \hbox{$c_1(\ga_J^*)\big|\bar{\cal U}_{{\cal T}_0(J')}\!=\!
c_1(\ga_{J'+\hat{1}}^*)$} 
under the decomposition~\e_ref{cart_split2},
the first claim of the lemma follows from~\e_ref{ag_lmm_e1}.\\
(2) Equation~\e_ref{ag_lmm_e1} implies that for any $\hat{1}\!\in\!J$,
\begin{equation}\label{ag_lmm_e2}
c_1(\ga_J^*)+\psi_{\hat{1}}=\!
\sum_{\eset\neq J'\varsubsetneq J-\{\hat{1}\}}\!\!\!\!\!\!
\hbox{PD}_{\bar{\cal U}_{{\cal T}_0}}
\bar{\cal U}_{{\cal T}_0(J')}.
\end{equation}
If ${\cal T}$, $i$, and $M_i$ are as in (2) of the lemma,
under the splitting~\e_ref{cart_split2},
\begin{equation}\label{ag_lmm_e3}
\hbox{PD}_{\bar{\cal U}_{\cal T}(\mu)}
\bar{\cal U}_{{\cal T}(M_i)}(\mu)\big|
\bar{\cal U}_{{\cal T}(M_i)}(\mu)=
-\psi_{\hat{1}}\!\times\!1-1\!\times\!\psi_{\hat{0}}.
\end{equation}
The second claim of the lemma follows from \e_ref{ag_lmm_e2},
applied with $J\!=\!\{\hat{1}\}\!+\!M_i$, 
and \e_ref{ag_lmm_e3}, since 
\hbox{$1\!\times\!\psi_{\hat{0}}\!=\!1\!\times\!c_1(L_i\bar{\cal T})$}.

\section{Comparison  of $n_d^{(1)}(\mu)$ and $n_{1,d}(\mu)$}
\label{degen_sect}

\subsection{Summary}
In this section, we prove

\begin{prp}
\label{g1degen_prp}
Suppose $n\!\ge\!2$, $d\!\ge\!1$, and $\mu\!=\!(\mu_1,\ldots,\mu_N)$ is 
an $N$-tuple of proper linear subspaces of~$\P$ in general position
such that \hbox{$\codim_{\Bbb{C}}\mu\!=\!d(n\!+\!1)\!-\!1$}.
Then 
$$n_d^{(1)}(\mu)=n_{1,d}(\mu).$$\\
\end{prp}

\noindent
Denote by $\ov{\frak M}_{1,1}$ the Deligne-Mumford moduli space
of stable genus-one curves with one marked point and 
by ${\frak M}_{1,1}$ the main stratum of~$\ov{\frak M}_{1,1}$,
i.e.~the complement of the point~$\i$ in~$\ov{\frak M}_{1,1}$.
The elements of ${\frak M}_{1,1}$ parameterize 
(equivalence classes of) smooth genus-one curves with one marked point.
The point $\i\!\in\!\ov{\frak M}_{1,1}$ corresponds 
to a sphere with one marked point and with two other points identified.\\

\noindent
Denote by $\ov{\frak M}\!=\!\ov{\frak M}_{1,N}\big(\P,d\big)$ 
the moduli space of 
stable \hbox{degree-$d$} maps
from $N$-pointed genus-one curves to~$\P$.
Let 
$$\ov{\frak M}(\mu)=
\big\{b\!\in\!\ov{\frak M}\!:\ev_l(b)\!\in\!\mu_l
~\forall l\!\in\![N]\big\}.$$
We denote by \hbox{$\pi\!:\ov{\frak M}\!\lra\!\ov{\frak M}_{1,1}$}
the forgetful functor sending each stable map
\hbox{$b\!=\![S,[N],I;x,(j,y),u]$}
to the one-marked curve $[S,y_1]$ and contracting 
all unstable components of~$(S,y_1)$.
The resulting complex curve is either a torus or 
a sphere with two points identified.
For any $\si\!\in\!\ov{\frak M}_{1,1}$, let
$$\ov{\frak M}_{\si}=\pi^{-1}(\si), \qquad
\ov{\frak M}_{\si}(\mu)=\ov{\frak M}_{\si}\cap \ov{\frak M}(\mu).$$
If the $j$-invariant~$\si$ is different from infinity,
i.e.~the stable curve~${\cal C}_{\si}$ corresponding to~$\si$ is smooth,
the cardinality of $\ov{\frak M}_{\si}(\mu)$
is $\big|\hbox{Aut}({\cal C}_{\si})\big|$ times
the number of genus-one degree-$d$ curves with $j$-invariant~$\si$
that pass through the constraints~$\mu$, i.e.
\begin{equation}
\label{g1degen_prp_e1}
\big|\ov{\frak M}_{\si}(\mu)\big|=2n_{1,d}(\mu).
\end{equation}
\\

\noindent
If $\big\{\si_k\big\}\!\subset\!{\frak M}_{1,1}$ converges to
$\i\!\in\!\ov{\frak M}_{1,1}$ and 
$b_k\!\in\!\ov{\frak M}_{\si_k}(\mu)$,
a subsequence of $\{b_k\}$ converges in $\ov{\frak M}$
to some \hbox{$b\!\in\!\ov{\frak M}_{\i}(\mu)$}.
It will be shown that $\Si_b$ is a sphere with two points identified;
see Lemma~\ref{dimen_count1} and Corollary~\ref{limits2}.
Conversely, for every 
$$b=(S,[N],\{\hat{0}\};,(\hat{0},y),u)\in\ov{\frak M}_{\i}(\mu)$$
such that $\Si_b$ is a sphere with two points identified
and for every $\si\!\in\!{\frak M}_{1,1}$ sufficiently close to~$\i$,
there exists a unique stable map $b(\si)\!\in\!\ov{\frak M}_{\si}(\mu)$
close to~$b$ in~$\ov{\frak M}$; see Lemma~\ref{dimen_count2}.
Since the number of stable maps 
$$b=(S,[N],\{\hat{0}\};,(\hat{0},y),u)\in\ov{\frak M}_{\i}(\mu)$$
such that $\Si_b$ is a sphere with two points identified
is $2n_d^{(1)}(\mu)$,
Proposition~\ref{g1degen_prp} follows from the two lemmas, 
the corollary,  and equation~\e_ref{g1degen_prp_e1}.

\subsection{Dimension Counts}

In this subsection, we show that if
$$[b]\!=\!\big[S,[N],I;x,(j,y),u\big]\in\ov{\frak M}_{\i}(\mu)$$
and $u_{\hat{0}}\!=\!u_b|S$ is not constant,
then $\Si_b\!=\!S$ is a sphere with two points identified;
see Lemma~\ref{dimen_count1}.
This lemma is proved by dimension counting.
We then observe that for each such stable map~$b$
and every $\si\!\in\!{\frak M}_{1,1}$ sufficiently close to~$\i$,
there exists a unique stable map $b(\si)\!\in\!\ov{\frak M}_{\si}(\mu)$
close to~$b$ in~$\ov{\frak M}$; see Lemma~\ref{dimen_count2}.

\begin{lmm}
\label{dimen_count1}
If $[b]\!=\!\big[S,[N],I;x,(j,y),u\big]\!\in\!\ov{\frak M}_{\i}(\mu)$ and
$u_{\hat{0}}\!=\!u_b|S$ is not constant,
then $\Si_b\!=\!S$ is a~sphere with two points identified.
\end{lmm}

\noindent
{\it Proof:} Suppose ${\cal T}\!=\!(S,[N],I;j,\under{d})$
is a simple bubble type such that $S$ is a circle of~$k$ spheres,
$d_{\hat{0}}\!\neq\!0$, and $\sum d_i\!=\!d$. 
Let ${\cal U}_{{\cal T},\under{d}_{\hat{0}}}$ denote the subspace
of ${\cal U}_{\cal T}$ such that the nonconstant restrictions of~$u_b$
to the components of $S$ have degrees $d_{\hat{0},1},\ldots,d_{\hat{0},k'}$
for all $b\!\in\!U_{{\cal T},\under{d}_{\hat{0}}}$.
We must have $\sum d_{\hat{0},l}\!=\!d_{\hat{0}}$.
Then, the dimension of $\dim{\cal U}_{{\cal T},\under{d}_{\hat{0}}}(\mu)$
is given~by
\begin{equation*}\begin{split}
\Big(\sum_{l=1}^{k'}\big(d_{\hat{0},l}(n\!+\!1)\!+\!n\!-\!1\big)-nk'
+\sum_{i\in\hat{I}}\big(d_i(n\!+\!1)\!+\!n\!-\!2\!-\!(n\!-\!1)\big)+
\big(N\!-\!(k\!-\!k')\big)\Big)
-\big(\codim_{\Bbb{C}}\mu\!+\!N)&\\
=1-|k|-|\hat{I}|&.
\end{split}\end{equation*}
Thus, 
${\cal U}_{\cal T}(\mu)\!=\!\eset$ unless $k\!=\!1$ and $\hat{I}\!=\!\eset$,
i.e.~$\Si_b\!=\!S$ is a sphere with two points identified.

\begin{lmm}
\label{dimen_count2}
For every $[b]\!=\!\big[S,[N],\{\hat{0}\};,(\hat{0},y),u\big]
\!\in\!\ov{\frak M}_{\i}(\mu)$
such that $S$ is a sphere with two points identified,
there exists neighborhood $U_b$ of $\i$ in $\ov{\frak M}_{1,1}$
and $W_b$ of $b$ in $\ov{\frak M}_{1,N}(\P,d)$ 
such~that 
$$\big|\ov{\frak M}_{\si}(\mu)\cap W_b\big|=1
\qquad\forall\si\!\in\!U_b\!-\!\{\i\}.$$
\end{lmm}

\noindent
{\it Proof:}
Since $d\!\ge\!1$,
\begin{equation}\label{dimen_count2e1}
H^1(S;u_b^*T\P)=(n\!+\!1)H^1\big(S;u_b^*{\cal O}(1_{\P})\big)=0,
\end{equation}
see Corollary~6.5 in~\cite{Z2} for example.
The lemma follows from~\e_ref{dimen_count2e1} by standard arguments.
A~purely analytic proof can be found in~\cite{RT}.

\subsection{A Property of Limits in $\ov{\frak M}_{1,N}(\P,d)$}

\noindent
Suppose $\big\{\si_k\big\}\!\subset\!{\frak M}_{1,1}$ converges to
$\i\!\in\!\ov{\frak M}_{1,1}$ and $b_k\!\in\!\ov{\frak M}_{\si_k}$
converges to
$$[b]\!=\!\big[S,[N],I;x,(j,y),u\big]\in\ov{\frak M}_{\i}$$
such that $u_b|S$ is constant.
In this subsection, we describe a condition such a limit~$b$
must satisfy; see Lemma~\ref{limits1}.
This lemma is the key part of Section~\ref{degen_sect}.
Its proof extends the argument of~\cite{P1} for the $n\!=\!2$ case 
and makes use of the explicit notation described in 
Subsection~\ref{notation}.
We conclude by  observing that no element of $\ov{\frak M}_{\i}(\mu)$
can satisfy this condition if the constraints~$\mu$
 are in general position.

\begin{lmm}
\label{limits1}
Suppose
$$[b]\!=\!\big[S,[N],I;x,(j,y),u\big]
\in \ov{\bigcup_{\si\in{\frak M}_{1,1}}\ov{\frak M}_{\si}}
\cap{\cal M}_{\cal T},$$
where ${\cal T}\!=\!(S,[N],I;j,\under{d})$ is a simple bubble type
such that $S$ is a circle of spheres and $d_{\hat{0}}\!=\!0$.
Then the dimension of the linear span of the set 
$\big\{du_h\big|_{\i}e_{\i}\!:h\!\in\!\chi({\cal T})\big\}$
is less than~$\big|\chi({\cal T})\big|$.
\end{lmm}

\noindent
{\it Proof:}
(1) By the algebraic geometry definition of stable-map convergence,
there exist\\
(i) a one-parameter family of curves
$\tilde{\ka}\!:\tilde{\cal F}\!\lra\!\De$
such that $\De$ is a neighborhood of~$0$ in~$\Bbb{C}$,
$\tilde{\cal F}$ is a~smooth space, $\tilde{\ka}^{-1}(0)\!=\!\Si_b$,
and $\Si_t\!\equiv\!\tilde{\ka}^{-1}(t)$ is a~smooth genus-one curve
for all $t\!\in\!\De^*\!\equiv\!\De\!-\!\{0\}$;\\
(ii) a holomorphic map $\tilde{u}\!:\tilde{\cal F}\!\lra\!\P$ such that
$\tilde{u}|\ka^{-1}(0)\!=\!u_b$.\\
This family $\tilde{\ka}\!:\tilde{\cal F}\!\lra\!\De$
can be obtained from another family of curves 
$\ka_{\hat{0}}\!\!:{\cal F}_{\hat{0}}\!\lra\!\De$
that satisfies~(i), except $\ka_{\hat{0}}^{-1}(0)\!=\!S$,
by a sequence of blowups at smooth points of the central fiber
as we now describe.
Choose an ordering~$\prec$ of the set~$I$ consistent 
with its partial ordering.
If $h\!\in\!I$, let
\begin{gather*}
I^h=\big\{i\!\in\!I\!: i\!\prec\!h\big\},
\quad i(h)=\max I^h\hbox{~~if~}h\!\in\!\hat{I},\quad
I^{(h)}=I^h\cup\{h\},\quad
M(h)=\big\{l\!\in\![N]: j_l\!\preceq\!h\big\},\\
b(h)=\big(S^2,M(h),I^{(h)};x|\hat{I}^{(h)},(j,y)|M(h),u|I^{(h)}\big).
\end{gather*}
Suppose $h\!\in\!\hat{I}$ and we have constructed 
a one-parameter family of curves 
$\ka_{i(h)}\!\!:{\cal F}_{i(h)}\!\lra\!\De$ that satisfies~(i),
except $\ka_{i(h)}^{-1}(0)\!=\!\Si_{b(i(h))}$.
Let ${\cal F}_h$ be the blowup of~${\cal F}_{i(h)}$ at 
the smooth point of $(\io_h,x_h)$ of~$\Si_{b(i(h))}$ and 
let $\ka_h\!\!:{\cal F}_h\!\lra\!\De$ be the induced projection map.
Choose coordinates $(t,w_h)$ near 
$(\io_h,x_h)\!\in\!{\cal F}_{i(h)}$ such that 
$d\ka_{i(h)}\frac{\partial}{\partial w_h}\!=\!0$,
i.e.~$w_h$ is a coordinate in 
$\ka_{i(h)}^{-1}(t)$ for $t\!\in\!\De$ sufficiently small.
We define coordinates $(t,z_h)$ on a neighborhood in~${\cal F}_h$
of the complement of the node of the new exceptional divisor~by
$$(t,z_h)\lra\big(t,w_h\!=\!tz_h,[1,z_h]\big).$$
For a good choice of the family 
$\ka_{\hat{0}}\!:{\cal F}_{\hat{0}}\!\lra\!\De$,
$\tilde{\cal F}\!=\!{\cal F}_{h^*}$ and $\tilde{\pi}\!=\!\pi_{h^*}$,
where $h^*$ is the largest element of $I$ with respect to
the ordering~$\prec$.\\
(2) Let $\psi\!\in\!H^0(S;w_S)$ be a nonzero differential,
i.e.~$\psi$ is a holomorphic $(1,0)$-form on the components of $S$,
which has simple poles at the singular points of~$S$
with residues that add up to zero at each node.
Then, for each $h\!\in\!H_{\hat{0}}{\cal T}$, there exists
$a_h\!\in\!\Bbb{C}^*$ such that
$$\psi|_{(0,w_h)}=a_h\big(1+o(1)\big)dw_h.$$
Thus, we can extend $\psi$ to a family of elements 
$\psi_t\!\in\!H^0(\Si_t;\om_{\Si_t})$ such that
\begin{equation}\label{limits1e1}
\psi|_{(t,w_h)}=a_h\big(1+o(1)\big)dw_h,
\quad\hbox{with}~~a_h\!\in\!\Bbb{C}^*.
\end{equation}
If $h\!\in\!\hat{I}$, let $|h|\!=\!\big|\{i\!\in\!I\!:i\!<\!h\}\big|$.
Denote by $\tilde{h}$ the element of $H_{\hat{0}}{\cal T}$
such that $h\!\in\!\bar{D}_{\tilde{h}}{\cal T}$.
By~\e_ref{limits1e1}, we have
\begin{equation}\label{limits1e2}
\psi|_{(t,z_h)}=t^{|h|}a_{\tilde{h}}\big(1+o(1_t)\big)dz_h,
\quad\hbox{with}~~a_{\tilde{h}}\!\in\!\Bbb{C}^*.
\end{equation}
(3) Let $H_1$ and $H_2$ be any two hyperplanes in~$\P$
that intersect the image of~$u_b$ transversally and 
miss the image of the nodes of~$\Si_b$.
Then for all~$t$ sufficiently small and $i\!=\!1,2$,
\begin{equation}\label{limits1e3}
u_t^{-1}(H_i)=\big\{z_{1,h_1}^{(i)}(t),\ldots,
                       z_{d,h_d}^{(i)}\big\}\subset\Si_t,
\quad\hbox{where}\quad h_j\!\in\!\hat{I},~~
z_{j,h_j}^{(i)}(t)=z_{j,h_j}^{(i)}(0)+o(1_t),
\end{equation}
$z_{j,h_j}^{(i)}(0)\!\in\!\Si_{b,h}$, and $u_t\!=\!\tilde{u}|\Si_t$.
Since $\sum z^{(1)}_{h_j}(t)$  and $\sum z^{(2)}_j(t)$
are linearly equivalent divisors in~$\Si_t$,
\begin{equation}\label{limits1e4}
\sum_{j=1}^{j=d}\int_{z^{(1)}_{j,h_j}(t)}^{z^{(2)}_{j,h_j}(t)}\!\psi_t=0
\qquad\forall t\!\in\!\De^*,
\end{equation}
where each line integral is taken inside of 
an appropriate coordinate chart~$(t,z_h)$.
Plugging \e_ref{limits1e2} and \e_ref{limits1e3}
into \e_ref{limits1e4} gives
\begin{equation}\label{limits1e5}
\sum_{j=1}^{j=d}t^{|h_j|}a_{\tilde{h}_j}
\big(z^{(2)}_{j,h_j}(0)-z^{(1)}_{j,h_j}(0)+o(1_t)\big)=0
\qquad\forall t\!\in\!\De^*.
\end{equation}
Let $k\!=\!\min\big\{|h|\!:h\!\in\!\chi({\cal T})\big\}$;
then $k\!=\!\min\big\{|h_j|\!:j\!\in\![d]\big\}$.
Thus, dividing equation~\e_ref{limits1e5} by~$t^k$
and then taking the limit as $t\!\lra\!0$, we conclude~that
\begin{equation}\label{limits1e6}
\sum_{|h_j|=k}\!a_{\tilde{h}_j}z^{(1)}_{j,h_j}(0)
=\sum_{|h_j|=k}\!a_{\tilde{h}_j}z^{(2)}_{j,h_j}(0).
\end{equation}
(4) Equality~\e_ref{limits1e6} holds for a dense subset of pairs~$(H_1,H_2)$.
The consequences of this fact can be interpreted as follows.
For each $h\!\in\!\hat{I}$, let $[u_h,v_h]$ be homogeneous
coordinates on $\Si_{b,h}$ such that \hbox{$z_h\!=\!v_h/u_h$}.
Each map $u_h$ corresponds to an $(n+1)$-tuple of
homogeneous polynomials
$$p_{h,i}=\sum_{l=0}^{l=d_h}p_{h,i;l}u^lv^{d-l},
\qquad i=0,\ldots,n,~~~ p_{h,i;l}\!\in\!\Bbb{C}.$$
Equality~\e_ref{limits1e6} implies that
there exists $K\!\in\!\Bbb{C}$ such~that
\begin{equation}\label{limits1e7}
\sum_{|h|=k,d_h\neq0}\!\!a_{\tilde{h}}
\frac{\sum_{i=0}^{i=n}c_ip_{h,i;d_h-1}}{\sum_{i=0}^{i=n}c_ip_{h,i;d_h}}=K
\qquad\forall\big[c_0,\ldots,c_n\big]\!\in\!\P.
\end{equation}
On the other hand, $u_{h_1}(\i)\!=\!u_{h_2}(\i)$ for all
$h_1,h_2\!\in\!\chi({\cal T})$.
Thus, for all  $h_1,h_2\!\in\!\chi({\cal T})$, there exists
$K_{h_1,h_2}\!\in\!\Bbb{C}^*\!-\!\{0\}$ such~that
$$\big(p_{h_1,0;d_{h_1}},\ldots,p_{h_1,n;d_{h_1}}\big)
=K_{h_1,h_2}\big(p_{h_2,0;d_{h_2}},\ldots,p_{h_2,n;d_{h_2}}\big).$$
It follows that \e_ref{limits1e7} is equivalent to
\begin{gather}
\sum_{i=0}^{i=n}\sum_{|h|=k,d_h\neq0}\!\!\!\!\tilde{a}_hp_{h,i;d_h-1}c_i
=K\sum_{i=0}^{i=n}p_{h_1,i;d_{h_1}}c_i\quad\forall c_i\!\in\!\Bbb{C}
\Lra\notag\\
\label{limits1e9}
\sum_{|h|=k,d_h\neq0}\!\!\!\!\tilde{a}_hp_{h,i;d_h-1}
=Kp_{h_1,i;d_{h_1}}, \qquad i\!=\!0,\ldots,n.
\end{gather}
where $h_1$ is a fixed element of the set 
$\big\{h\!\in\!\hat{I}\!:|h|\!=\!k,d_h\!\neq\!0\big\}$ and
$\tilde{a}_h\!\in\!\Bbb{C}^*$.
It is straightforward to deduce from~\e_ref{limits1e9} that
$$\sum_{|h|=k,d_h\neq0}\!\!\!\tilde{a}_hdu_h\big|_{\i}e_{\i}=0.$$
The lemma is now proved, since 
$\big\{h\!\in\!\hat{I}\!:|h|\!=\!k,d_h\!\neq\!0\big\}
\!\subset\!\chi({\cal T})$.

\begin{crl}
\label{limits2}
Suppose
$$[b]\!=\!\big[S,[N],I;x,(j,y),u\big]
\in \ov{\bigcup_{\si\in{\frak M}_{1,1}}\ov{\frak M}_{\si}}
\cap\ov{\frak M}_{\i}(\mu).$$
Then $u_b|S$ is not constant.
\end{crl}

\noindent
{\it Proof:} Suppose $u_b|S$ is constant. Let
\begin{gather*}
\tilde{I}\!=\!\big\{i\!\in\!I\!:\chi_{\cal T}i\!\neq\!0\big\}
\!\subset\!\hat{I},~~ 
M_0\!=\!\!\bigcup_{i\in I-\tilde{I}}\!\!\!M_i{\cal T},~~
\tilde{x}\!=\!x|\hat{\tilde{I}},~~
(\tilde{j},\tilde{y})\!=\!(j,y)\big|\big([N]\!-\!M_0\big),~~
\tilde{d}\!=\!d|\tilde{I},~~
\tilde{u}\!=\!u|\tilde{I};\\
\tilde{\cal T}=\big(S^2,[N]\!-\!M_0,\tilde{I};\tilde{j},\tilde{d}\big),
\quad
\tilde{b}=\big(S^2,[N]\!-\!M_0,\tilde{I};
\tilde{x},(\tilde{j},\tilde{y}),\tilde{u}\big).
\end{gather*}
Then, $\tilde{\cal T}$ is a bubble type such that $\sum\tilde{d}_i\!=\!d$
and $\tilde{d}_i\!>\!0$ for all $i\!\in\!\tilde{I}\!-\!\hat{\tilde{I}}$.
The latter property implies that 
$\chi(\tilde{\cal T})\!=\!\tilde{I}\!-\!\hat{\tilde{I}}$.
Furthermore, $\tilde{b}\!\in\!{\cal U}_{\tilde{\cal T}}(\mu)$.
By Lemma~\ref{limits1}, the linear map
$$\al_{|\chi(\tilde{\cal T})|,|M_0|}\!:
\bigoplus_{i\in\chi(\tilde{\cal T})}\!\!\!L_i\tilde{\cal T}\lra\ev^*T\P,
\quad
\al_{|\chi(\tilde{\cal T})|,|M_0|}(\ups)
=\sum_{i\in\chi(\tilde{\cal T})}\!\!{\cal D}_{\tilde{\cal T},i}\ups_i,$$
does not have full rank at~$\tilde{b}$.
However, this is impossible by Lemma~\ref{str_lmm}.

\end{document}